\documentclass{amsart}
\usepackage{amsmath,amssymb, amsbsy}
\usepackage{color,psfrag}
\usepackage[dvips]{graphicx}
\usepackage{enumerate}
\newcommand{\R}{{\mathbb R}}
\newcommand{\N}{{\mathbb N}}

\renewcommand{\H}{{\mathcal H}}
\newcommand{\SN}{{\mathbb S}^{N-1}}

\newcommand{\e }{\varepsilon}

\def\ds{\displaystyle}

\renewcommand{\ge }{\geqslant}
\renewcommand{\geq }{\geqslant}
\renewcommand{\le }{\leqslant}
\renewcommand{\leq }{\leqslant}

\newenvironment{pf}{\noindent{\sc Proof}.\enspace}{\hfill\qed\medskip}
\newenvironment{pfn}[1]{\noindent{\bf Proof of
    {#1}.\enspace}}{\hfill\qed\medskip}
\newtheorem{Theorem}{Theorem}[section]
\newtheorem{Corollary}[Theorem]{Corollary}
\newtheorem{Lemma}[Theorem]{Lemma}
\newtheorem{Proposition}[Theorem]{Proposition}
\theoremstyle{definition}

\newtheorem{remark}[Theorem]{Remark}
\begin{document}

\title[Borderline Hardy potentials]{On semilinear elliptic equations\\
  with borderline Hardy potentials}

\author[Veronica Felli]{Veronica Felli}
\address{\hbox{\parbox{5.7in}{\medskip\noindent{Universit\`a degli
        Studi di Milano-Bicocca,\\
        Dipartimento di Ma\-t\-ema\-ti\-ca e Applicazioni, \\
        Via Cozzi
        53, 20125 Milano, Italy. \\[3pt]
        \em{E-mail address: }{\tt veronica.felli@unimib.it.}}}}}
\author[Alberto Ferrero]{Alberto Ferrero}
\address{\hbox{\parbox{5.7in}{\medskip\noindent{Universit\`a degli Studi del Piemonte Orientale,\\
        Viale Teresa Michel 11, 15121 Alessandria, Italy. \\[3pt]
        \em{E-mail addresses: }{\tt
          alberto.ferrero@mfn.unipmn.it}.}}}}  

\date{September 21, 2012.}

\keywords{Hardy's inequality, singular  elliptic operators,
   asymptotic behavior of solutions.}

\subjclass[2010]{35J75, 35B40, 35B45.}

\thanks{V. Felli was
partially supported by the PRIN2009 grant ``Critical Point Theory and
Perturbative Methods for Nonlinear Differential
Equations''. A. Ferrero was partially supported by the PRIN2008 grant
  ``Aspetti geometrici delle equazioni alle derivate parziali e
  questioni connesse''.}

\begin{abstract}
   \noindent
   In this paper we study the asymptotic behavior of solutions to an
   elliptic equation near the singularity of an inverse square
   potential with a coefficient related to the best constant for the
   Hardy inequality. Due to the presence of a borderline Hardy
   potential, a proper variational setting has to be introduced in
   order to provide a weak formulation of the equation.  An
   Almgren-type monotonicity formula is used to determine the exact
   asymptotic behavior of solutions.
\end{abstract}

\maketitle

\section{Introduction}\label{intro}
On a domain $\Omega\subseteq \R^N$, $N\ge 3$, containing the
origin, let us consider the following problem
\begin{equation} \label{eq:u}
-\Delta u-\left(\frac{N-2}2\right)^{\!\!2} \frac{u}{|x|^2}=h(x) u+f(x,u), \quad \text{in } \Omega,
\end{equation}
where $h$ is possibly singular at the origin but negligible with
respect to the Hardy potential and $f$ is a nonlinearity subcritical
with respect to the critical Sobolev exponent. Looking at equation
\eqref{eq:u}, one may observe that the best constant for the
classical Hardy inequality appears in front of the inverse square
potential; this can be considered as a borderline situation for
several points of view, from the variational setting to the
existence and qualitative behavior of solutions. Recent papers
were devoted to equations and differential inequalities involving
elliptic operators with inverse square potentials in the
borderline situation, see
\cite{bre-dupai-tesei,cirstea,dupaigne,fall-musina,garcia-peral-primo,tesei,
vazquez_zuazua}.

In \cite{fall-musina} the authors study necessary conditions for the
existence of nonnegative distributional solutions of the differential
inequality
\begin{equation*}
  -\Delta u-\left(\frac{N-2}2\right)^{\!\!2} \frac{u}{|x|^2}\ge \alpha \frac{u}{|x|^2\log^2 |x|} \quad
  \text{in } \mathcal D'(B_R\setminus \{0\}),
\end{equation*}
where $B_R$ denotes the ball of radius $R$ centered at the origin. The logarithmic term appearing in the above inequality is related to an improved version of the Hardy inequality, see for example \cite{adim-chaud-rama,barbi-filippa-terti,chaud,ghous-morad}.

In \cite{cirstea}, the author studies existence of positive distributional solutions of the nonlinear elliptic equation
\begin{equation*}
-\Delta u-\lambda \frac{u}{|x|^2}+b(x)h(u)=0 \quad \text{in }
\Omega\setminus \{0\} ,
\end{equation*}
satisfying some prescribed asymptotic behaviors at the origin, where $\Omega\subset\R^N$, $N\ge 3$ is a domain containing the origin and $\lambda\in (-\infty,(N-2)^2/4]$. These prescribed asymptotic behaviors are related to the following fundamental solutions
$$
\Phi_\lambda^{+}(x)=|x|^{-\frac{N-2}2-\sqrt{\frac{(N-2)^2}4-\lambda}}   \, , \qquad
\Phi_\lambda^-(x)=|x|^{-\frac{N-2}2+\sqrt{\frac{(N-2)^2}4-\lambda}} \, ,
$$
$$
\Psi^{+}(x)=|x|^{-\frac{N-2}2}\log(1/|x|) \, , \qquad \Psi^-(x)=|x|^{-\frac{N-2}2}
$$
of
$$
-\Delta u-\lambda \frac{u}{|x|^2}=0 \qquad \text{in } \R^N\setminus\{0\}
$$
respectively in the cases $\lambda\in(-\infty,(N-2)^2/4)$ and $ \lambda=(N-2)^2/4$.

Similar results were obtained in \cite{bra-chia-cirstea-trombetti} for
equations with elliptic operators in divergence form.

In \cite{dupaigne} the author studies a singular elliptic Dirichlet
problem with a power type nonlinearity and a forcing term:
\begin{equation} \label{eq:dupaigne}
\begin{cases}
\ds{-\Delta u-\frac{c}{|x|^2}\, u=u^p+tf } & \text{in } \Omega,\\
u>0 & \text{in } \Omega, \\
u=0 & \text{on } \partial\Omega,
\end{cases}
\end{equation}
where $\Omega\subset\R^N$, $N\ge 3$, is a domain containing the
origin, $p>1$, $t>0$, $f$ a smooth, bounded, nonnegative function, and
$c\in (0,(N-2)^2/4]$. In \cite{dupaigne}, the author provides a
classification of different kind of solutions of problem
\eqref{eq:dupaigne}, both of distributional and variational type.  In
the present paper, we are going to introduce an analogous terminology
for solutions to \eqref{eq:u}, a classification of which will be
provided as a byproduct of our main result, see section \ref{s:main}
for details.

We also mention that, in \cite{bre-dupai-tesei} the authors study existence and
nonexistence of solutions of the equation in \eqref{eq:dupaigne}
with $t=0$.

In the spirit of \cite{dupaigne}, in the present paper we concentrate
our attention on local solutions to \eqref{eq:u} belonging to a
suitable functional space related to the borderline case of the Hardy
inequality. To this purpose, in section \ref{s:main} we introduce the
Hilbert space $\mathcal H(\omega)$ defined as completion of $C^\infty_c(\overline\omega\setminus\{0\})$ 
with respect to a scalar product related to
the Hardy potential appearing in \eqref{eq:u} (see \eqref{eq:H(omega)}). Here $\omega$
represents a bounded domain with $\partial\omega\in C^1$.

The purpose of this paper is to classify the possible asymptotic
behaviors of solutions to \eqref{eq:u} near the singularity of the
Hardy potential. Some results in this direction were obtained in
\cite{FFT4,FFT,FFT2,FFT3} for different kinds of problems: in
\cite{FFT, FFT3} Schr\"odinger equations with electromagnetic
potentials, in \cite{FFT2} Schr\"odinger equations with inverse
square many-particle potentials; finally in \cite{FFT4} the
authors study the asymptotic behavior of solutions of a singular
elliptic equation near a corner of the boundary. For other results
concerning elliptic equations with singular inverse square
potentials see also \cite{FMT2,FG,Jan,Han-P.,terracini96}.

The results of the present paper are closely related to
the ones obtained in \cite{FFT,FFT3}. If we drop the magnetic part of the
electromagnetic potential, the equation studied in \cite{FFT, FFT3}
becomes
\begin{equation} \label{eq:H^1}
-\Delta u-\frac{a(x/|x|)}{|x|^2} u=h(x)u+f(x,u) \quad \text{in }
\Omega,
\end{equation}
where $a\in L^\infty(\SN)$. In \cite{FFT,FFT3} the quadratic form
associated to the linear
operator $-\Delta-\frac{a(x/|x|)}{|x|^2}$ is assumed to satisfy a coercivity
type condition. More precisely, it is required that the first
eigenvalue $\mu_1(0,a)$ of the spherical operator
$-\Delta_{\SN}-a(\theta)$ satisfies
$\mu_1(0,a)>-\left(\frac{N-2}2\right)^{\!2}$.

Due to this coercivity, it was quite natural in that setting
looking for $H^1$-solutions to \eqref{eq:H^1}, i.e. functions $u\in
H^1(\Omega)$ satisfying \eqref{eq:H^1} in a variational sense, whereas, 
 in the borderline situation considered in the present paper,
it is reasonable to replace the classical $H^1$ Sobolev space with
the above mentioned $\mathcal H$ space.

In the proof of our main result (Theorem \ref{t:asymptotic} below) we
perform an Almgren-type monotonicity procedure (see \cite{almgren,GL})
and provide a characterization of the leading term in the asymptotic
expansion by means of a Cauchy's integral type representation formula.

As an application of the main result, we also prove an a priori
estimate and a unique continuation principle for solutions to
\eqref{eq:u}; see \cite{Kurata} for questions related to unique
continuation principles for elliptic equations with singular
potentials.

This paper is organized as follows. In Section \ref{s:main} we
introduce the assumptions of the main result and explain in details
what we mean by a $\mathcal H$-solution of \eqref{eq:u}. In Section
\ref{s:var-form} we describe the main properties of the space
$\mathcal H$, while in Section \ref{sec:an-equiv-probl} we reformulate
\eqref{eq:u} in cylindrical variables, introducing an auxiliary
equation in a cylinder of $\R^{N+1}$. In Section \ref{s:Almgren} we
study the Almgren-type function associated to the problem, which is
combined in Section \ref{sec:blow-up-argument} with a blow-up argument
to characterize
the leading term in the asymptotic expansion of solutions of
\eqref{eq:u} near the origin, thus proving the main theorem.

\bigskip

\section{Assumptions and main results}\label{s:main}
We first introduce the assumptions on the
potential $h$ and the nonlinearity $f$.
We assume that $h$ satisfies
\begin{equation} \label{eq:h}\tag{$\bf H$}
h\in L^{\infty}_{\rm loc}(\Omega\setminus\{0\}), \quad |h(x)|\leq
C_h|x|^{-2+\e}  \quad \mbox{in $\Omega\setminus\{0\}$ for some
$C_h>0$ and $\e>0$}.
\end{equation}
It is not restrictive to assume that $\e\in(0,2)$.
Let  $f$ satisfy
\begin{equation}\label{F}\tag{$\bf F$}
\left\{\!\!
\begin{array}{l}
f\in C^0(\Omega\times \R),\quad F\in C^1(\Omega\times \R),
\quad s\mapsto f(x,s)\in C^1(\R)\quad\text{for a.e. }x\in\Omega, \\[5pt]
|f(x,s)s|+|f'_s(x,s)s^2|+|\nabla_x F(x,s)\cdot x|\leq C_f(|s|^2+|s|^{p})
\quad\,\text{for a.e. $x\in\Omega$ and all $s\in\R$},
\end{array}
\right.
\end{equation}
 where $F(x,s)=\int_0^s f(x,t)\,dt$,
$2<p<2^*=\frac{2N}{N-2}$, $C_f>0$ is a
constant independent of $x\in\Omega$ and $s\in\R$, $\nabla_x F$
denotes the gradient of $F$ with respect to the $x$ variable, and
$f'_s(x,s)=\frac{\partial f}{\partial s}(x,s)$.

In order to state the main result of this paper, a suitable
variational formulation for solutions of \eqref{eq:u} has to be
introduced (see also  \cite{dupaigne,garcia-peral-primo,
vazquez_zuazua}).
For any bounded domain $\omega\subset \R^N$ containing the origin
and satisfying $\partial \omega\in C^1$, let us define
$\H(\omega)$ as the completion of the space
$C^\infty_c(\overline\omega\setminus\{0\})$ with respect to the
scalar product
\begin{align}\label{eq:H(omega)}
(u,v)_{\H(\omega)}:=&\int_{\omega}
  \nabla u(x)\cdot \nabla v(x)\, dx-\left(\frac{N-2}{2}\right)^2
  \int_{\omega} \frac{u(x)v(x)}{|x|^2}\, dx \\
  & \notag +\int_{\omega} u(x)v(x)\, dx +\frac{N-2}2
  \int_{\partial\omega} \frac{u(x)v(x)}{|x|^2} \, (x\cdot\nu(x))\, dS,
  \quad u,v\in C^\infty_c(\overline\omega\setminus\{0\}).
\end{align}
The form in \eqref{eq:H(omega)} is actually a scalar product on
$C^\infty_c(\overline\omega\setminus\{0\})$ as 
detailed in Section \ref{s:var-form}.

For any domain $\Omega\subseteq \R^N$ satisfying
$0\in \Omega$ (with $\partial\Omega$ not necessarily in $C^1$), we
define  the
space $\H_{{\rm loc}}(\Omega)$ as the space of functions $u\in
H^1_{{\rm loc}}(\Omega\setminus\{0\})$ such that $u_{|\omega}\in
\H(\omega)$ for any domain $\omega\Subset \Omega$ with
$\partial\omega\in C^1$.

We are ready to provide a rigorous definition for solutions to
\eqref{eq:u}. Let $h$, $f$ satisfy respectively \eqref{eq:h} and
\eqref{F}: by a solution of \eqref{eq:u} we mean a function
$u\in\H_{{\rm loc}}(\Omega)$ such that
\begin{equation} \label{eq:variational} \int_{\Omega} \nabla
  u(x)\cdot\nabla \varphi(x)\, dx-\left(\frac{N-2}2\right)^{\!\!2}
  \int_{\Omega} \frac{u(x)}{|x|^2}\, \varphi(x)\, dx =\int_{\Omega}
  \big(h(x) u(x)+f(x,u(x))\big)\varphi(x)\, dx
\end{equation}
for any $\varphi\in C^\infty_c (\Omega\setminus\{0\})$.
We observe that every term in the above identity is well-defined in view of Proposition \ref{p:H1loc} and Proposition
\ref{p:embedding}.

The above notion of solution corresponds to the notion of
$\mathcal H(\Omega)$-solution introduced in \cite[Section
6]{dupaigne}, as we will prove in Proposition \ref{p:2.6}. In other
words, if $u\in \mathcal H_{{\rm loc}}(\Omega)$ is an $\mathcal
H$-solution of \eqref{eq:u}, then for any $\omega\Subset \Omega$
with $\partial\omega\in C^1$ we have
\begin{equation*}
(u,v)_{\H(\omega)}=\int_{\omega} (h(x)+1)u(x)v(x)\, dx+
\int_{\omega} f(x,u(x))v(x)\, dx \quad \text{for any } v\in
\mathcal H_0(\omega),
\end{equation*}
where $\mathcal H_0(\omega)$ is the closure in $\mathcal H(\omega)$ of
the space $C^\infty_c(\omega\setminus\{0\})$.

In \cite{dupaigne}, the following notion of
strong solution is also discussed. By a strong solution to \eqref{eq:u} we mean a
function $u\in C^2(\Omega\setminus\{0\})$ which solves
\eqref{eq:u} in the classical sense and satisfies the following
pointwise estimate: for any $R>0$ there exists a constant
$C=C(N,h,f,u,\Omega,R)$ depending only on $N,h,f,u,\Omega,R$ but
independent of $x$ such that
\begin{equation*}
|u(x)|\leq C |x|^{-\frac{N-2}2} \quad \text{for any } x\in
(\Omega\cap B_R)\setminus \{0\} .
\end{equation*}
Before giving the statement of our main result, we recall that
the eigenvalues of the Laplace Beltrami operator $-\Delta _{\mathbb{S}^{N-1}}$
are given by
\begin{equation*}
\lambda _{\ell }=(N-2+\ell )\ell ,\quad \ell =0,1,2,\dots ,
\end{equation*}%
having the $\ell $-th eigenvalue $\lambda _{\ell }$ multiplicity
\begin{equation*}
m_{\ell }=\frac{(N-3+\ell )!(N+2\ell -2)}{\ell !(N-2)!},
\end{equation*}%
and the eigenfunctions coincide with the usual spherical harmonics. For
every $\ell \geq 0$, let $\{Y_{\ell ,m}\}_{m=1}^{m_{\ell }}$ be
a
$L^{2}(\mathbb{S}^{N-1})$-orthonormal basis of the
eigenspace of $-\Delta _{\mathbb{S}^{N-1}}$ associated to $\lambda _{\ell }$
with $Y_{\ell ,m}$ being spherical harmonics of degree $\ell $.

\begin{Theorem} \label{t:asymptotic}
Let $N\ge 3$ and assume \eqref{eq:h}, \eqref{F}. Let $u\in
\mathcal H_{{\rm loc}}(\Omega)$ be a nontrivial $\mathcal
H$-solution to~(\ref{eq:u}). Then there exist $\ell_0\in \N$ and
$\beta_{\ell_0,1},\dots,\beta_{\ell_0,m_{\ell_0}}\in \R$ such that
$(\beta_{\ell_0,1},\dots,\beta_{\ell_0,m_{\ell_0}})\neq(0,\dots,0)$
and, for
any $\alpha\in (0,1)$,
\begin{equation*}
r^{\frac{N-2}2-\sqrt{\lambda_{\ell_0}}} \,u(r\theta)\to
\sum_{m=1}^{m_{\ell_0}} \beta_{\ell_0,m} Y_{\ell_0,m}(\theta) \quad
\text{in } C^{1,\alpha}(\SN) \quad \text{as } r\to 0^+
\end{equation*}
and
\begin{equation*}
r^{\frac{N}2-\sqrt{\lambda_{\ell_0}}} \nabla u(r\theta)\to
\sum_{m=1}^{m_{\ell_0}} \beta_{\ell_0,m}
\left[\left(-\tfrac{N-2}2+\sqrt{\lambda_{\ell_0}}\right)
Y_{\ell_0,m}(\theta)\theta+\nabla_{\SN}
Y_{\ell_0,m}(\theta)\right]
\end{equation*}
in $C^{0,\alpha}(\SN)$ as $r\to 0^+$.
Moreover,  the coefficients
$\beta_{\ell_0,1},\dots,\beta_{\ell_0,m_{\ell_0}}$ admit the following
representation
\begin{align*}
\beta_{\ell_0,m}=\int_{\SN} \!\! \left[\!
\frac{u(R\theta)}{R^{\widetilde \gamma}}+ \!\!\int_0^R
\frac{h(s\theta)u(s\theta)+f(s\theta,u(s\theta))}{2\widetilde
\gamma+N-2}\left(\!s^{-\widetilde\gamma+1}-\frac{s^{\widetilde\gamma+N-1}}{R^{2\widetilde\gamma+N-2}}\right)
 \! ds \right] Y_{\ell_0,m}(\theta)\, dS(\theta)
\end{align*}
for any $R>0$ such that $B_R:=\{x\in \R^N:|x|<R\}\subset \Omega$,
where $\widetilde\gamma:=-\frac{N-2}2+\sqrt{\lambda_{\ell_0}}$.
\end{Theorem}

\noindent As a consequence of Theorem \ref{t:asymptotic}, the following
pointwise estimates hold true.
\begin{Corollary} \label{c:pointwise} Let $N\ge 3$ and assume
  \eqref{eq:h}, \eqref{F}. If $u\in \mathcal H_{{\rm loc}}(\Omega)$
  is a nontrivial $\mathcal H$-solution to~(\ref{eq:u}), then there
  exists $\ell_0\in \N$ such that
\begin{equation*}
|u(x)|=O\Big(|x|^{-\frac{N-2}2+\sqrt{\lambda_{\ell_0}}}\Big)
\quad\text{and}\quad
|\nabla u(x)|=O\Big(|x|^{-\frac{N}2+\sqrt{\lambda_{\ell_0}}}\Big)
\quad\text{as }|x|\to 0.
\end{equation*}
\end{Corollary}
\noindent The following result follows immediately from Corollary \ref{c:pointwise}.

\begin{Corollary}\label{c:strong}
  Let $N\ge 3$ and assume \eqref{eq:h}, \eqref{F}. Let $u\in \mathcal
  H_{{\rm loc}}(\Omega)$ be a nontrivial $\mathcal H$-solution
  to~(\ref{eq:u}). Then
\begin{enumerate}[(i)]
\item $|u(x)|=O\big(|x|^{-\frac{N-2}2}\big)$ as $|x|\to 0$. 
\item If $u$ changes sign in a neighborhood of $0$, then $u\in
  H^1_{\rm loc}(\Omega)$.
\end{enumerate}
\end{Corollary}

From classical elliptic regularity theory and Corollary \ref{c:strong},
it follows that if $u$ is an $\mathcal H$-solution
 to (\ref{eq:u}) and $h,f$ are smooth outside $0$, then $u$ is a \emph{strong solution} in the sense of \cite{dupaigne}.

As another byproduct of Theorem \ref{t:asymptotic}, we also have
the following version of the \emph{Strong Unique Continuation Principle} for an
elliptic equation with a singular coefficient.
\begin{Corollary} Let $N\ge 3$ and assume \eqref{eq:h}, \eqref{F}. Let $u\in
\mathcal H_{{\rm loc}}(\Omega)$ be a $\mathcal H$\!-solution of
(\ref{eq:u}). If $u(x)=O(|x|^k)$ as $|x|\to 0^+$ for any $k\in \N$,
then $u\equiv 0$ in $\Omega$.
\end{Corollary}

\medskip
\noindent
{\bf Notation. } \par
\begin{itemize}
\item[-] For all $r>0$, $B_r$ denotes the open ball $\{x\in
\R^N:|x|<r\}$ in $\R^N$ with center at $0$ and radius $r$.
\item[-] $C^\infty_c(A)$ denotes the space of
$C^\infty(A)$-functions whose support is compact in $A$. \item[-]
For any open set $\Omega\subseteq \R^N$, $\mathcal D'(\Omega)$
denotes the space of distributions on $\Omega$.
 \item[-]  $dS$ denotes the volume
element on the spheres $\partial B_r$, $r>0$.
 \item[-] For any $N\ge 1$ we
put $\omega_{N-1}:=\int_{\SN} dS$.
\end{itemize}

\section{On $\mathcal H$-solutions to \eqref{eq:u}} \label{s:var-form}
In this section we describe the main properties of
the space $\mathcal H$ and of $\mathcal H$-solutions to \eqref{eq:u}.
In the sequel, $\omega$ denotes a bounded domain in $\R^N$
satisfying $\partial\omega\in C^1$ and $0\in\omega$.
In order to reformulate \eqref{eq:u} in cylindrical variables, we let 
\begin{equation*}
\Phi:\R^N\setminus\{0\}\to \mathcal C:=\R\times \SN\subset \R^{N+1}
\end{equation*}
be the diffeomorphism \emph{(Emden-Fowler transformation)} defined as
\begin{equation*}
\Phi(x):=\left(-\log|x|,\frac{x}{|x|}\right) \quad \text{for any } x\in \R^N\setminus\{0\},
\end{equation*}
see \cite{CW2}, and let $\mathcal
C_\omega:=\Phi(\omega\setminus\{0\})\subseteq \mathcal C$. Let us
introduce the linear operator
\begin{align}\label{eq:Tu}
  &T:C^\infty_c (\overline\omega\setminus\{0\})\to
  C^\infty_c(\overline{\mathcal C}_\omega),\\
  \notag& Tu(t,\theta):=e^{-\frac{N-2}2t}u(e^{-t}\theta), \quad \text{
    for any } (t,\theta)\in \mathcal C,\,u\in C^\infty_c
  (\overline\omega\setminus\{0\}) .
\end{align}
Clearly $T$ is an isomorphism between vector spaces.  Let us denote by
$\mu$ the standard volume measure on the cylinder $\mathcal C$, by
$\nabla_{\mathcal C}$ the gradient associated with the standard Riemannian
metric of $\mathcal C$, and by $(t,\theta)$ the generic element of
$\mathcal C$.

We observe that $(\cdot,\cdot)_{\H(\omega)}$ as defined in
\eqref{eq:H(omega)} is actually a scalar product on
$C^\infty_c(\overline\omega\setminus\{0\})$ since the following
identities hold for any $u\in
C^\infty_c(\overline\omega\setminus\{0\})$, see \cite{CW1,CW2,wz}:
\begin{align}\label{eq:id1}
  &\int_{\omega} |\nabla u(x)|^2 dx-\left(\frac{N-2}{2}\right)^2 \int_{\omega} \frac{u^2(x)}{|x|^2}\, dx
 \\
&\notag\qquad\qquad  +\frac{N-2}2 \int_{\partial\omega} \frac{u^2(x)}{|x|^2} \,
  (x\cdot\nu(x))\, dS=\int_{\mathcal C_\omega} |\nabla_{\mathcal C}(Tu)|^2 d\mu\ge 0,\\
&\label{eq:id2}
\int_{\omega} u^2(x) \, dx=\int_{\mathcal C_\omega} e^{-2t} (Tu)^2 d\mu.
\end{align}
We also define, for any $\omega$ as above, the weighted Sobolev space
$H_\mu(\mathcal C_\omega)$ as the completion of $C^\infty_{\rm
  c}(\overline{\mathcal C}_\omega)$ with respect to the norm
\begin{equation*}
\|w\|_{H_\mu(\mathcal C_\omega)}:=\left(\int_{\mathcal C_\omega}
|\nabla_{\mathcal C} w|^2 d\mu +\int_{\mathcal C_\omega} e^{-2t}
w^2 d\mu\right)^{1/2}.
\end{equation*}
By density and continuity it is possible to extend $T$ as a linear
and continuous operator from $\H(\omega)$ to $H_\mu(\mathcal
C_\omega)$. In this way $T:\H(\omega)\to H_\mu(\mathcal C_\omega)$
becomes an isometric isomorphism.

The following proposition relates $\H(\omega)$ with the classical
Sobolev space $H^1(\omega)$.
\begin{Proposition} \label{p:H1}
Let $\omega\subset \R^N$ be a bounded domain satisfying $0\in
\omega$ and $\partial\omega\in C^1$. Then $H^1(\omega)\subset
\H(\omega)$ with continuous embedding.
\end{Proposition}
We omit the proof of Proposition \ref{p:H1} which can be easily obtained by 
classical density  arguments. We observe the
inclusion $H^1(\omega)\subset
\H(\omega)$ is actually strict, since the function
$|x|^{-\frac{N-2}{2}}$ belongs to $\H(\omega)$ but not to $H^1(\omega)$.


\begin{Proposition} \label{p:H1loc}
Let $\omega\subset \R^N$ be a bounded domain satisfying $0\in
\omega$ and $\partial\omega\in C^1$. Then $\H(\omega)\subseteq
H^1_{{\rm loc}}(\overline\omega\setminus\{0\})\cap L^2(\omega)$
where by $H^1_{{\rm loc}}(\overline\omega\setminus\{0\})$ we mean
the space of functions which belong to $H^1(A)$ for any open set
$A$ satisfying $\overline A\subseteq
\overline\omega\setminus\{0\}$.
\end{Proposition}

\begin{pf}
Let $u\in \H(\omega)$ and let $\{u_n\}\subset C^\infty_c(\overline\omega\setminus\{0\})$ be a sequence such that
 $u_n\to u$ in $\H(\omega)$. Then, for any $n$, $Tu_n$ belongs to $H_\mu(\mathcal C_\omega)$. By
\eqref{eq:Tu}, \eqref{eq:id1}, \eqref{eq:id2} and direct calculations,
for any open set $A$ such that $\overline A\subseteq
\overline\omega\setminus\{0\}$, we have that, denoting $V_{n,m}=T(u_n-u_m)$,
\begin{align} \label{eq:un-um}
& \|u_n-u_m\|^2_{H^1(A)}=\int_{\mathcal C_A}
\bigg[|\nabla_{\mathcal C} V_{n,m}|^2
+\big(\tfrac{N-2}2\big)^{\!2} V_{n,m}^2\bigg] d\mu
+\int_{\mathcal C_A} \bigg[\tfrac{N-2}2 \,\partial_t (V_{n,m}^2)+e^{-2t}V_{n,m}^2 \bigg] d\mu \\
\notag & \le \int_{\mathcal C_A} \bigg[|\nabla_{\mathcal C}
V_{n,m}|^2
+\big(\tfrac{N-2}2\big)^{\!2} V_{n,m}^2\bigg] d\mu +\int_{\mathcal C_A} \bigg[\tfrac{N-2}2 |\partial_t
V_{n,m}|^2+\tfrac{N-2}2 V_{n,m}^2
+e^{-2t}V_{n,m}^2 \bigg] d\mu \\
\notag & \le K_A \bigg(\int_{\mathcal C_A}
|\nabla_{\mathcal C} V_{n,m}|^2 d\mu+\int_{\mathcal C_A} e^{-2t} V_{n,m}^2 d\mu\bigg) \\
\notag & \le K_A \bigg(\int_{\mathcal C_\omega} |\nabla_{\mathcal
C} V_{n,m}|^2 d\mu+\int_{\mathcal C_\omega} e^{-2t}
V_{n,m}^2 d\mu\bigg) =K_A \|u_n-u_m\|^2_{\H(\omega)}
\end{align}
where $K_A:=\max\left\{ \frac{N}2, 1+\frac{N(N-2)}4
\sup_{(t,\theta)\in\mathcal C_A} e^{2t}\right\}$.
Then $\{u_n\}$ is a Cauchy sequence in $H^1(A)$ for any
$A$ as above and hence $u\in\H(\omega)$ may be seen as a function
in $H^1_{{\rm loc}}(\omega\setminus\{0\})$.  Moreover, by
\eqref{eq:H(omega)}, \eqref{eq:id1} and \eqref{eq:id2} it is clear
that
$$
\|u_n-u_m\|_{L^2(\omega)}\leq \|u_n-u_m\|_{\H(\omega)}
$$
and hence $\{u_n\}$ is a Cauchy sequence in $L^2(\omega)$. This
completes the proof of the proposition.~\end{pf}

\noindent In \cite[Extension 4.3]{BV} the following Poincar\'e-Sobolev type
inequality was proved.

\begin{Proposition} \label{p:Poincare}
\ \cite[Extension 4.3]{BV} Let $\omega\subset \R^N$ be a bounded
domain satisfying $0\in \omega$ and let $1\le q<\frac{2N}{N-2}$.
Then there exists a constant $C(\omega,q)$ such that
\begin{equation} \label{eq:Poincare}
\left(\int_{\omega} |u(x)|^q \, dx\right)^{\!\!2/q} \leq C(\omega,q)
\left[\int_{\omega} |\nabla u(x)|^2 dx -\left(\frac{N-2}2\right)^{\!\!2}
\int_\omega \frac{u^2(x)}{|x|^2} \, dx \right] 
\end{equation}
for any $u\in C^\infty_c(\omega\setminus\{0\})$. 
\end{Proposition}

Let us consider the space $\H_0(\omega)$ defined in Section \ref{s:main} as the closure in $\H(\omega)$ of
$C^\infty_c(\omega\setminus\{0\})$, see 
also \cite[Section 6]{dupaigne}. If 
we define the scalar product
\begin{equation}\label{eq:H_0}
  (u,v)_{\H_0(\omega)}:=\int_{\omega} \nabla u(x)\cdot \nabla v(x)\,
  dx-
\left(\frac{N-2}{2}\right)^{\!\!2} \int_{\omega} \frac{u(x)v(x)}{|x|^2}\, dx,
  \quad \text{for any } u,v\in C^\infty_c(\omega\setminus\{0\}),
\end{equation}
then by (\ref{eq:Poincare}) with $q=2$ we deduce that the norms
$\|\cdot\|_{\H(\omega)}$ and $\|\cdot\|_{\H_0(\omega)}$ are
equivalent on $C^\infty_c(\omega\setminus\{0\})$. Hence $\H_0(\omega)$ may be endowed with
the equivalent scalar product obtained by density, extending the
scalar product $(\cdot,\cdot)_{\H_0(\omega)}$ defined in
\eqref{eq:H_0} to the whole $\H_0(\omega)\times \H_0(\omega)$. 

By Proposition \ref{p:Poincare} and the definition of
$\H_0(\omega)$, the following Sobolev type embedding follows.

\begin{Proposition}\label{p:Poincare2}
Let $\omega\subset \R^N$ be a bounded domain satisfying $0\in
\omega$ and let $1\le q<\frac{2N}{N-2}$. Then $\H_0(\omega)\subset
L^q(\omega)$ with continuous embedding.
\end{Proposition}

Actually the continuous embedding
$\H(\omega)\subset L^q(\omega)$, $1\le q<2N/(N-2)$, also holds true as
shown Proposition \ref{p:embedding}.

\begin{Proposition} \label{p:embedding} Let $\omega\subset \R^N$ be a bounded domain satisfying $0\in
\omega$ and let $1\le q<\frac{2N}{N-2}$. Then $\H(\omega)\subset
L^q(\omega)$ with continuous embedding.
\end{Proposition}

\begin{pf} Let $u\in \H(\omega)$. Then by Proposition
\ref{p:H1loc} we deduce that $u\in L^q(\omega\setminus \overline
B_\delta)$ for any $\delta>0$ such that $\overline B_\delta\subset
\omega$. Moreover, arguing as in \eqref{eq:un-um}, we infer that there exists
a constant $C(N,q,\delta)$ depending only on $N,q,\delta,\omega$ such
that
\begin{equation} \label{eq:embedding-1}
\|u\|_{L^q(\omega\setminus B_\delta)}\leq C(N,q,\delta)
\|u\|_{\H(\omega)} \, .
\end{equation}
Let us prove that for some fixed $\delta>0$ chosen as above $u\in
L^q(B_\delta)$. To this purpose let $\eta\in C^\infty_c(\omega)$
be a radial function such that $0\leq \eta\leq 1$ and $\eta\equiv 1$ in $\overline
B_\delta$. Let $\{u_n\}\subset
C^\infty_c(\overline\omega\setminus\{0\})$ be a sequence such that
$u_n\to u$ in $\H(\omega)$. Then by \eqref{eq:id1},
\eqref{eq:id2}, \eqref{eq:Poincare} we obtain
\begin{align} \label{eq:Cauchy}
 \bigg(\int_\omega |&(u_n-u_m)\eta|^q dx\bigg)^{\!\!\frac2q} \!\!\le C(\omega,q)\! \left[\int_\omega
|\nabla((u_n-u_m)\eta)|^2 dx-\big(\tfrac{N-2}2\big)^{\!2}\!\!
\int_\omega \frac{|(u_n-u_m)\eta|^2}{|x|^2}\, dx\right]
\\
\notag & =C(\omega,q) \int_{\mathcal C_\omega}
|\nabla_{\mathcal C} (T((u_n-u_m)\eta))|^2 d\mu  \\
\notag &
\leq 2 \int_{\mathcal C_\omega} \eta^2(e^{-t}\theta)|\nabla_{\mathcal C} (T(u_n-u_m))|^2 d\mu
+2\int_{\mathcal C_\omega} |T(u_n-u_m)|^2 |\nabla_{\mathcal C}(\eta(e^{-t}\theta))|^2 d\mu
\\
\notag & \leq 2 \int_{\mathcal C_\omega} |\nabla_{\mathcal C}
(T(u_n-u_m))|^2 d\mu +2 \|\nabla \eta\|_{L^\infty(\omega)}^2
\int_{\mathcal C_\omega} e^{-2t} |T(u_n-u_m)|^2 d\mu \\
\notag & \leq 2(1+\|\nabla \eta\|_{L^\infty(\omega)}^2)
\|u_n-u_m\|^2_{\H(\omega)} \, .
\end{align}
This shows that $\{u_n\eta\}$ is a Cauchy sequence in
$L^q(\omega)$. Since $u_n\eta\to
u\eta$ pointwise then $u\eta\in L^q(\omega)$. In particular $u\in L^q(B_\delta)$. Moreover proceeding as in
\eqref{eq:Cauchy} we also have that
\begin{align} \label{eq:embedding-2}
 \|u\|_{L^q(B_\delta)}& \leq
\|\eta u\|_{L^q(\omega)} =\lim_{n\to +\infty} \|\eta u_n\|_{L^q(\omega)} \\
\notag &
\leq \lim_{n\to +\infty} [2(1+\|\nabla
\eta\|_{L^\infty(\omega)}^2)]^{1/2} \|u_n\|_{\H(\omega)}
=[2(1+\|\nabla \eta\|_{L^\infty(\omega)}^2)]^{1/2} \|u\|_{\H(\omega)} \, .
\end{align}
Combining \eqref{eq:embedding-1} and \eqref{eq:embedding-2} we
conclude that $\H(\omega)\subset
L^q(\omega)$ with continuous embedding.
\end{pf}

\noindent From Propositions \ref{p:H1loc} and
 \ref{p:embedding} we infer that, if $u\in \H_{{\rm
loc}}(\Omega)$, then $u\in H^1_{{\rm
loc}}(\Omega\setminus\{0\})\cap L^q_{{\rm loc}}(\Omega)$ for all
$1\leq q<2N/(N-2)$.

\begin{remark}\label{rem:reg}
 From Proposition \ref{p:H1loc}, we have that if $u\in \H_{{\rm loc}}(\Omega)$
is a solution to (\ref{eq:u}) in the sense of~(\ref{eq:variational}), then $u$ is a weak $H^1$-solution in
$\Omega\setminus\{0\}$. Hence, classical Brezis-Kato \cite{BK} estimates, bootstrap,
  and elliptic regularity theory, imply that
$u\in H^2_{\rm loc}(\Omega\setminus \{0\})\cap
C^{1,\alpha}_{\rm loc}(\Omega\setminus\{0\})$ for any $\alpha\in
(0,1)$.
\end{remark}

\noindent From \eqref{eq:variational} we deduce the following characterizations of solutions to
\eqref{eq:u}.
\begin{Proposition}\label{p:2.6} Let $h$ satisfy \eqref{eq:h},  $f$ satisfy \eqref{F} and 
$u\in \H_{{\rm loc}}(\Omega)$ be a solution to (\ref{eq:u}) in the sense of
(\ref{eq:variational}). Then $u$ solves (\ref{eq:u}) in the sense
of distributions in $\Omega$, i.e.
\begin{equation} \label{eq:distributional}
-\int_{\Omega} u(x)\Delta\varphi(x)\,
dx-\left(\frac{N-2}2\right)^{\!\!2} \int_{\Omega} \frac{u(x)}{|x|^2}\,
\varphi(x)\, dx =\int_{\Omega}
 \big(h(x) u(x)+f(x,u(x))\big)
\varphi(x)\, dx
\end{equation}
for any $\varphi\in C^\infty_c (\Omega)$.
Moreover, for any bounded domain $\omega$ with $\partial \omega\in
C^1$ and $\overline \omega\subset \Omega$, we have that
\begin{equation} \label{eq:variationalB}
(u,v)_{\H(\omega)}=\int_{\omega} (h(x)+1)u(x)v(x)\, dx+
\int_{\omega} f(x,u(x))v(x)\, dx
\end{equation}
for all $v\in \H_0(\omega)$ and
\begin{equation} \label{eq:variational2}
  (u,v)_{\H(\omega)}=\int_{\omega} (h+1)uv\, dx+
  \int_{\omega} f(x,u)v\, dx
  +\int_{\partial \omega} \frac{\partial u}{\partial \nu}v \, dS
  +\frac{N-2}2\int_{\partial \omega} \frac{uv}{|x|^2} \, (x\cdot \nu)
  \, dS
\end{equation}
for all $v\in \H(\omega)$.
\end{Proposition}

\begin{pf} Let $\varphi\in C^\infty_c(\Omega)$ and let $\omega$ be an
  open domain with smooth boundary satisfying $\overline\omega\subset
  \Omega$ and $\text{supp}\, \varphi\subset \omega$. Since $u\in
  \H_{{\rm loc}}(\Omega)$ then $u\in \H(\omega)$ and there exists a
  sequence $\{u_n\}\subset C^\infty_c(\overline \omega\setminus\{0\})$
  such that $u_n\to u$ in $\H(\omega)$. Since $\varphi\in
  H^1_0(\omega)$, Proposition \ref{p:H1} implies that there exists a sequence
  $\{\varphi_m\}\subset C^\infty_c(\omega\setminus\{0\})$ such that
  $\varphi_m\to \varphi$ in $H^1(\omega)$ and in $\H(\omega)$ as
  $m\to+\infty$. Hence we have
\begin{align} \label{eq:scalar_product}
& (u,\varphi)_{\H(\omega)}=\lim_{n\to +\infty}\left(\lim_{m\to+\infty} (u_n,\varphi_m) _{\H(\omega)}\right)\\
\notag & =\lim_{n\to +\infty}\bigg\{\lim_{m\to +\infty} \bigg[\int_\omega
\nabla u_n\cdot \nabla\varphi_m \, dx
-\left(\frac{N-2}2\right)^{\!\!2} \int_{\omega}
\frac{u_n\varphi_m}{|x|^2}\, dx +\int_{\omega}
u_n\varphi_m\, dx\bigg]\bigg\}\\
\notag & =\lim_{n\to +\infty}\left[\int_\omega \nabla u_n\cdot
\nabla\varphi \, dx -\left(\frac{N-2}2\right)^{\!\!2} \int_{\omega}
\frac{u_n\varphi}{|x|^2}\, dx +\int_{\omega}
u_n\varphi\, dx\right] \\
\notag & =\lim_{n\to +\infty}\left[-\int_\omega u_n\Delta \varphi
\, dx -\left(\frac{N-2}2\right)^{\!\!2} \int_{\omega}
\frac{u_n\varphi}{|x|^2}\, dx +\int_{\omega}
u_n\varphi\, dx\right].
\end{align}
By Proposition \ref{p:embedding} we also have that $u_n\to u$ in
$L^q(\omega)$ for any $1\leq q<2^*$. By H\"older inequality with $\frac{N}{N-2}<p<\frac{2N}{N-2}$ we
have that $\frac{u\varphi}{|x|^2}\in L^1(\omega)$ and
\begin{align*}
\left|\int_{\omega} \frac{u_n(x)\varphi(x)}{|x|^2}\, dx-
\int_{\omega} \frac{u(x)\varphi(x)}{|x|^2}\, dx\right|
\leq \|\varphi\|_{L^\infty(\omega)} \left(\int_\omega |u_n(x)-u(x)|^q dx\right)^{1/q}
\left(\int_{\omega} |x|^{-\frac{2q}{q-1}} dx\right)^{\frac{q-1}q}
\end{align*}
and hence passing to the limit in \eqref{eq:scalar_product} we obtain
\begin{equation} \label{eq:identity1}
(u,\varphi)_{\H(\omega)}=-\int_\omega u(x)\Delta \varphi(x) \, dx
-\left(\frac{N-2}2\right)^2 \int_{\omega}
\frac{u(x)\varphi(x)}{|x|^2}\, dx +\int_{\omega} u(x)\varphi(x)\,
dx                                \, .
\end{equation}
On the other hand, by the convergence $u_n\to u$ in $H^1_{{\rm
loc}}(\overline\omega\setminus\{0\})$, see Proposition
\ref{p:H1loc}, we obtain
\begin{align} \label{eq:scalar_product_2}
(u&,\varphi)_{\H(\omega)}=\lim_{m\to +\infty}\left(\lim_{n\to +\infty}(u_n,\varphi_m)_{\H(\omega)}\right) \\
\notag & =\lim_{m\to +\infty}\bigg\{\lim_{n\to +\infty}
\bigg[\int_\omega \nabla u_n\cdot \nabla\varphi_m \, dx
-\left(\frac{N-2}2\right)^{\!\!2} \int_{\omega}
\frac{u_n\varphi_m}{|x|^2}\, dx +\int_{\omega}
u_n\varphi_m\, dx\bigg]\bigg\} \\
\notag & =\lim_{m\to +\infty} \left[\int_\omega \nabla u\cdot
\nabla\varphi_m \, dx -\left(\frac{N-2}2\right)^{\!\!2} \int_{\omega}
\frac{u\varphi_m}{|x|^2}\, dx +\int_{\omega}
u\varphi_m\, dx\right]  \\
\notag & =\lim_{m\to +\infty}\bigg[ \int_\omega
(h+1)u\varphi_m\, dx+
\int_{\omega} f(x,u)\varphi_m\, dx\bigg] 
= \int_\omega
(h+1)u\varphi\, dx+
\int_{\omega} f(x,u)\varphi\, dx
\end{align}
where the last identity follows from
assumptions \eqref{eq:h} and \eqref{F} and the fact that $\varphi_m\to
\varphi$ in $L^q(\omega)$ for any $1\leq q<\frac{2N}{N-2}$. Combining \eqref{eq:identity1}
and \eqref{eq:scalar_product_2} obtain \eqref{eq:distributional}. 

The proof of \eqref{eq:variationalB} follows by the following density
argument: let $\{v_m\}\subset C^\infty_c(\omega \setminus \{0\})$ such
that $v_m\to v$ in $\H(\omega)$. Now it is enough to pass to the limit
as $m\to +\infty$ in \eqref{eq:scalar_product_2} with $v_m$ in place
of $\varphi$.

It remains to prove \eqref{eq:variational2}.  By elliptic regularity
estimates $u\in C^1(\overline\omega\setminus \{0\})$ (see Remark \ref{rem:reg}) and hence the
normal derivative of $u$ on $\partial\omega$ is continuous.  Let
$v\in \H(\omega)$ and let $\{v_m\}\subset
C^\infty_c(\overline\omega\setminus\{0\})$ be such that $v_m\to v $ in
$\H(\omega)$. Therefore we are allowed to integrate by parts to obtain
\begin{align*}
(u&,v_m)_{\H(\omega)}=\!\int_{\omega} \nabla u\cdot \nabla v_m \, dx-\left(\frac{N-2}2\right)^{\!\!2}\!\!
\int_{\omega} \frac{u}{|x|^2} \, v_m\, dx +\!\int_{\omega} uv_m\, dx+
\frac{N-2}2\! \int_{\partial \omega} \frac{uv_m}{|x|^2} (x\cdot \nu)\, dS \\
=& \int_{\omega} -(\Delta u)v_m \, dx+\int_{\partial \omega} \frac{\partial u}{\partial \nu}
v_m \, dS-\left(\frac{N-2}2\right)^{\!\!2}
\int_{\omega} \frac{u}{|x|^2} \, v_m\, dx \\
&\quad\quad+\int_{\omega} uv_m\, dx+
\frac{N-2}2 \int_{\partial \omega} \frac{uv_m}{|x|^2} (x\cdot \nu)\, dS \\
 =& \int_\omega huv_m \, dx
+
\int_{\omega} f(x,u)v_m\, dx
+\int_{\partial \omega} \frac{\partial u}{\partial \nu} v_m \, dS+\int_{\omega} uv_m\, dx+
\frac{N-2}2 \int_{\partial \omega} \frac{uv_m}{|x|^2} (x\cdot \nu)\, dS \, .
\end{align*}
The proof of \eqref{eq:variational2} follows passing to the limit as $m\to+\infty$.
\end{pf}

\section{An equivalent problem on the cylinder $\mathcal C$}\label{sec:an-equiv-probl}
 Reformulation of \eqref{eq:variational} in cylindrical
variables yields the following characterizations of solutions to
\eqref{eq:u}.

\begin{Proposition} Let $h$ satisfy \eqref{eq:h}, $f$ satisfy \eqref{F}, and  $u\in \H_{{\rm loc}}(\Omega)$ be a
solution to (\ref{eq:u}) in the sense of (\ref{eq:variational}). 
If $\omega$ is a bounded domain with $\partial \omega\in C^1$
and $\overline\omega\subset \Omega$, then the function $v:=Tu\in
H_\mu(\mathcal C_\omega)$ is a weak solution of the equation
\begin{equation} \label{eq:equation_in_C}
-\Delta_{\mathcal C} v(t,\theta)=e^{-2t}\widetilde
h(t,\theta)v(t,\theta)+e^{-2t}\widetilde f(t,\theta,v(t,\theta)),
 \quad \text{in } \mathcal C_\omega,
\end{equation}
where $\Delta_{\mathcal C}$ denotes the Laplace-Beltrami operator on $\mathcal C$ and
\begin{align*}
\widetilde h(t,\theta):=h(e^{-t}\theta),
\quad
  \widetilde
  f(t,\theta,s):=e^{-\frac{N-2}{2}t}f\Big(e^{-t}\theta,e^{\frac{N-2}{2}t}s\Big),
\quad\text{for any $(t,\theta)\in \mathcal C_\omega$},
\end{align*}
in the sense that 
\begin{equation} \label{eq:variational_H0mu}
\int_{\mathcal C_\omega} \nabla_{\mathcal C} v \cdot
\nabla_{\mathcal C} w \, d\mu =\int_{\mathcal C_\omega} e^{-2t}
\big( \widetilde h v+\widetilde f(t,\theta,v)\big) w \, d\mu,
\quad \text{for every } w\in H_{\mu,0}(\mathcal
C_\omega):=T(\H_0(\omega)).
\end{equation}
Moreover
\begin{equation} \label{eq:variational_Hmu}
\int_{\mathcal C_\omega} \nabla_{\mathcal C} v \cdot \nabla_{\mathcal C} w \, d\mu
=\int_{\partial \mathcal C_\omega} (\nabla_{\mathcal C}v \cdot \nu_{\partial \mathcal C_\omega}) \, w
\, dS+\int_{\mathcal C_\omega} e^{-2t} \big( \widetilde h v+\widetilde f(t,\theta,v)\big)
w \, d\mu,
\end{equation}
for every $w\in H_{\mu}(\mathcal C_\omega)$, where
$\nu_{\partial\mathcal C_\omega}$ denotes the exterior normal
vector to $\partial\mathcal C_\omega$ on $\mathcal C$.
\end{Proposition}

\noindent The following corollary is an immediate consequence of \eqref{eq:variational_Hmu}.

\begin{Corollary} Let $h$ satisfy \eqref{eq:h}, let $f$ satisfy \eqref{F}, and
let $u\in \H_{{\rm loc}}(\Omega)$ be a solution to (\ref{eq:u}) in
the sense of (\ref{eq:variational}). For any $t\in\R$, let
\begin{equation*}
\mathcal C_t:=\{(s,\theta)\in \mathcal C:s>t, \theta\in \SN\}, \quad 
\Gamma_t:=\{(t,\theta)\in\mathcal C:\theta\in\SN\}.
 \end{equation*}
Then for any $t$ such that $\overline{\mathcal C}_t\subset
\mathcal C_\Omega$, the function $v:=Tu\in H_\mu(\mathcal C_t)$
satisfies
\begin{multline}\label{eq:variational_Hmu_bis}
  \int_{\mathcal C_t} \nabla_{\mathcal C} v \cdot \nabla_{\mathcal C}
  w \, d\mu =-\int_{\Gamma_t} \frac{\partial v}{\partial s}
  \, w \, dS
  +\int_{\mathcal C_t} e^{-2s} \big( \widetilde h(s,\theta)
  v(s,\theta)+\widetilde f(s,\theta,v(s,\theta))\big) w(s,\theta) \,
  d\mu
\end{multline}
for any $w\in H_{\mu}(\mathcal C_t)$.
\end{Corollary}

\noindent In order to study solutions to
\eqref{eq:equation_in_C}, the properties of space
$H_\mu$ have to be investigated. The next results go in this direction.

\begin{Lemma}\label{l:equiv_norm}
For every $t\in\R$,   $H_{\mu}(\mathcal C_t)\hookrightarrow
L^2(\Gamma_t)$ with compact embedding. Furthermore,
\begin{equation}\label{eq:1}
v\mapsto \left(\int_{\mathcal C_t}
|\nabla_{\mathcal C} v|^2 d\mu +\int_{\Gamma_t}v^2dS\right)^{\!\!1/2}
\end{equation}
is an equivalent norm in $H_{\mu}(\mathcal C_t)$; more precisely,
there exists a constant $C>0$ such that, for all $t\in\R$ and $v\in
H_{\mu}(\mathcal C_t)$,
\begin{align}\label{eq:3}
\frac1C\bigg(
\int_{\mathcal C_t}|\nabla_{\mathcal C}v|^2\,d\mu+e^{2t}\int_{\mathcal
      C_t}e^{-2s}v^2\,d\mu\bigg)&\leq
\int_{\mathcal C_t}
|\nabla_{\mathcal C} v|^2 d\mu +\int_{\Gamma_t}v^2dS\\
\notag&\leq C\bigg(
\int_{\mathcal C_t}|\nabla_{\mathcal C}v|^2\,d\mu+e^{2t}\int_{\mathcal
      C_t}e^{-2s}v^2\,d\mu\bigg).
  \end{align}
\end{Lemma}
\begin{pf}
  The embedding  $H_{\mu}(\mathcal C_t)\hookrightarrow
L^2(\Gamma_t)$  and its compactness are just a consequence of the fact
  that   $T:\H(\omega)\to H_\mu(\mathcal C_\omega)$ is an isometric
  isomorphism combined with Proposition \ref{p:H1loc} and compactness
  of classical Sobolev trace embeddings. To show that the quadratic
  form in \eqref{eq:1} is an equivalent norm in $H_{\mu}(\mathcal
  C_t)$, we notice that, for all $v\in C^\infty_{\rm
    c}(\overline{\mathcal C}_t)$, integration by parts yields
  \begin{align}\label{eq:2}
    \int_{\mathcal
      C_t}e^{-2s}v^2(s,\theta)\,d&\mu(s,\theta)=\int_{\SN}\bigg(\int_t^{+\infty}e^{-2s}v^2(s,\theta)ds\bigg)dS(\theta)\\
\notag&=\int_{\SN}\bigg(\bigg[-\frac12
e^{-2s}v^2(s,\theta)\bigg]_{s=t}^{s=+\infty}+\int_t^{+\infty}e^{-2s}
\frac{dv}{ds}(s,\theta)v(s,\theta)ds\bigg)dS(\theta)\\
\notag&=\frac12 e^{-2t}\int_{\Gamma_t}v^2dS+\int_{\mathcal C_t}e^{-2s}
v\frac{dv}{ds}\,d\mu,
  \end{align}
which implies
\begin{align*}
  \int_{\mathcal
      C_t}e^{-2s}v^2\,d\mu\leq \frac12 e^{-2t}\int_{\Gamma_t}v^2dS+
\frac12 \int_{\mathcal
      C_t}e^{-2s}v^2\,d\mu+\frac{e^{-2t}}2
\int_{\mathcal C_t}|\nabla_{\mathcal C}v|^2\,d\mu
\end{align*}
and hence
\begin{align*}
  \int_{\mathcal
      C_t}e^{-2s}v^2\,d\mu\leq  e^{-2t}\bigg(\int_{\Gamma_t}v^2dS+
\int_{\mathcal C_t}|\nabla_{\mathcal C}v|^2\,d\mu\bigg).
\end{align*}
On the other hand, \eqref{eq:2} also implies
\begin{align*}
  e^{-2t}\int_{\Gamma_t}v^2dS\leq 3 \int_{\mathcal
      C_t}e^{-2s}v^2\,d\mu+e^{-2t}\int_{\mathcal C_t}|\nabla_{\mathcal C}v|^2\,d\mu.
\end{align*}
The conclusion then follows by density.
\end{pf}

\noindent The following lemma provides a Hardy type inequality
with boundary terms.

\begin{Lemma}\label{l:hardy_bound}
  For every $\sigma>0$ and $t\in\R$, $H_{\mu}(\mathcal C_t)\subset
  L^2(\mathcal C_t,e^{-\sigma s}d\mu)$. Furthermore,  for every $\sigma>0$
  there exists $\widetilde C_\sigma>0$ such that, for all $t\in\R$ and $v\in H_{\mu}(\mathcal C_t)$,
$$
\int_{\mathcal C_t}e^{-\sigma s}v^2(s,\theta)d\mu(s,\theta)\leq \widetilde C_\sigma e^{-\sigma
  t}\bigg(\int_{\mathcal C_t}| \nabla_{\mathcal C} v
|^2d\mu+\int_{\Gamma_t}v^2 dS\bigg).
$$
\end{Lemma}
\begin{pf}
For all $v\in C^\infty_{\rm
    c}(\overline{\mathcal C}_t)$, integration by parts yields
  \begin{align*}
    \int_{\mathcal
      C_t}e^{-\sigma
      s}v^2(s,\theta)\,d\mu(s,\theta)&=\int_{\SN}\bigg(\int_t^{+\infty}e^{-\sigma
      s}v^2(s,\theta)ds\bigg)dS(\theta)\\
\notag&=\int_{\SN}\bigg(\bigg[-\frac1\sigma
e^{-\sigma
  s}v^2(s,\theta)\bigg]_{s=t}^{s=+\infty}+\frac2\sigma\int_t^{+\infty}e^{-\sigma
  s}
\frac{dv}{ds}(s,\theta)v(s,\theta)ds\bigg)dS(\theta)\\
\notag&=\frac1\sigma e^{-\sigma t}\int_{\Gamma_t}v^2dS+\frac2\sigma\int_{\mathcal
  C_t}e^{-\sigma s}
v\frac{dv}{ds}\,d\mu,
  \end{align*}
which implies
\begin{align*}
  \int_{\mathcal
      C_t}e^{-\sigma s}v^2\,d\mu\leq \frac1\sigma e^{-\sigma t}\int_{\Gamma_t}v^2dS+
\frac12 \int_{\mathcal
  C_t}e^{-\sigma s}v^2\,d\mu+\frac{2e^{-\sigma t}}{\sigma^2}
\int_{\mathcal C_t}|\nabla_{\mathcal C}v|^2\,d\mu
\end{align*}
and hence
\begin{align*}
  \int_{\mathcal
      C_t}e^{-\sigma s}v^2\,d\mu\leq  e^{-\sigma t  }\bigg(\frac2\sigma\int_{\Gamma_t}v^2dS+\frac4{\sigma^2}
\int_{\mathcal C_t}|\nabla_{\mathcal C}v|^2\,d\mu\bigg).
\end{align*}
The conclusion thereby follows with $\widetilde C_\sigma=\max\{2/\sigma,4/\sigma^2\}$.
\end{pf}

\noindent The following Hardy-Sobolev type inequality holds.

\begin{Lemma}\label{l:sob_cil}
  For every $q\in\big[1,\frac{2N}{N-2}\big)$, there exists $C_{N,q}>0$
  such that, for all $t\in\R$ and $v\in H_{\mu}(\mathcal C_t)$,
$$
\bigg(\int_{\mathcal C_t}e^{(-N+\frac{N-2}{2}q)s}|v(s,\theta)|^q
d\mu(s,\theta)\bigg)^{\!\! 2/q}\leq C_{N,q}\, e^{(-\frac
  {2N}q+N-2)t}\bigg(\int_{\mathcal C_t}| \nabla_{\mathcal C} v
|^2d\mu+\int_{\Gamma_t}v^2 dS\bigg).
$$
\end{Lemma}
\begin{pf}
  From Proposition \ref{p:embedding}, there exists $c_{N,q}>0$ such
  that
$$
\bigg(\int_{B_1}|u(x)|^q dx\bigg)^{\!\!1/q}\leq
c_{N,q}\|u\|_{\mathcal
  H(B_1)}
$$
for all $u\in \mathcal H(B_1)$. Performing the change of variable
$v(s,\theta)=Tu(s-t,\theta)$ in the above inequality for all $t\in
\R$ and taking into account  \eqref{eq:3}, we obtain the stated
inequality with $C_{N,q}=c_{N,q}^2C$.\end{pf}

\section{The Almgren frequency function}\label{s:Almgren}

In this section, our purpose would be to construct an appropriate Almgren-type
frequency function for the solution to problem
\eqref{eq:variational}. Since for a general function
$u\in\H(\omega)$, the norm $\|u\|_{\H(\omega)}$ cannot be
expressed in an integral form, we prefer to look for an
Algrem-type function associated with the function $v:=Tu$.

In a domain $\Omega\subset \R^N$, let $u\in \H_{{\rm loc}}(\Omega)$ be a solution
of (\ref{eq:variational}). Let $R>0$ be such that $\overline B_R\subset \Omega$.
According with \cite{almgren,GL} (see also
\cite{FFT4,FFT,FFT2,FFT3}), for $t>-\log R$,  we define the functions
\begin{equation} \label{eq:D(t)}
D(t):=\int_{\mathcal C_t} |\nabla_{\mathcal C} v|^2 d\mu-\int_{\mathcal C_t} e^{-2s} \widetilde
hv^2 d\mu-\int_{\mathcal C_t} e^{-2s} \widetilde f(s,\theta,v)v\,d\mu
\end{equation}
and
\begin{equation} \label{H(t)}
H(t):=\int_{\Gamma_t} v^2 dS,
\end{equation}
where  $v:=Tu$ and $T$ is defined in \eqref{eq:Tu}.

\begin{Lemma}\label{welld} Let $\Omega\subset \R^N$ a domain a let $u\in \H_{{\rm loc}}(\Omega)$ be a solution
of (\ref{eq:variational}),  $u\not\equiv 0$, with
  $h$ satisfying \eqref{eq:h} and $f$ satisfying \eqref{F}. Let $H=H(t)$ be the function defined in
  (\ref{H(t)}). Then there exists $\bar t>0$ such that $H(t)>0$
  for any $t>\bar t$.
\end{Lemma}
\begin{pf}
Let us argue by contradiction and assume that there exists
$t_n\to+\infty$ such that $H(t_n)=0$; in particular $v=0$ on
$\Gamma_{t_n}$ and $v\in H_{\mu,0}(\mathcal C_{t_n})$. From
\eqref{eq:variational_Hmu_bis}, \eqref{eq:h}, \eqref{F}, and Lemmas
\ref{l:hardy_bound} and \ref{l:sob_cil}
\begin{align*}
0&=\int_{\mathcal C_{t_n}} |\nabla_{\mathcal C} v|^2 d\mu-\int_{\mathcal C_{t_n}} e^{-2s} \widetilde
hv^2 d\mu-\int_{\mathcal C_{t_n}} e^{-2s} \widetilde f(s,\theta,v)v\,d\mu\\
&
\geq \int_{\mathcal C_{t_n}} |\nabla_{\mathcal C} v|^2 d\mu-C_h
\int_{\mathcal C_{t_n}} e^{-\e s}v^2d\mu-C_f
\int_{\mathcal C_{t_n}} e^{-2 s}v^2d\mu-C_f
\int_{\mathcal C_{t_n}}
e^{\big(-N+\frac{N-2}{2}p\big)s}|v|^pd\mu\\
&\geq \Big(1-C_h\widetilde C_\e e^{-\e t_n}-C_f\widetilde C_2  e^{-2
  t_n}-C_f
C_{N,p}^{p/2}\, e^{(-N+\frac
  {N-2}2p)t_n}\big({\textstyle{\int_{\mathcal C_{-\log R}}| \nabla_{\mathcal C} v
|^2d\mu}}\big)^{\!\frac{p-2}{2}}\Big) \int_{\mathcal C_{t_n}}
|\nabla_{\mathcal C} v|^2 d\mu\\
&=(1+o(1)) \int_{\mathcal C_{t_n}}
|\nabla_{\mathcal C} v|^2 d\mu
  \end{align*}
which implies that $v\equiv 0$ in $\mathcal C_{t_n}$ for $n$
sufficiently large. Hence $u\equiv 0$ in a neighborhood of the origin
and, by classical
unique continuation principles for second order elliptic equations
with locally bounded coefficients (see e.g.  \cite{wolff}) we
conclude that $u=0$ a.e. in $\Omega$, a contradiction.
\end{pf}

\noindent By virtue of Lemma \ref{welld}, the \emph{Almgren-type
frequency
  function}
\begin{equation}\label{N(t)}
{\mathcal N}(t)=\frac{D(t)}{H(t)}
\end{equation}
is well defined in $(\bar t,+\infty)$.

In order to obtain a suitable representation for the derivative of $D$ we need the following Pohozaev-type
identity.
\begin{Proposition}\label{p:Pohozaev} 
  Let $h$ satisfy \eqref{eq:h}, $f$ satisfy \eqref{F}, and 
  $u\in \H_{{\rm loc}}(\Omega)$ be a solution to
  (\ref{eq:variational}). Let $R>0$ be such that $\overline B_R\subset
  \Omega$. Then for every $t\in (-\log R,+\infty)$ the function
  $v:=Tu\in H_\mu(\mathcal C_{-\log R})$ satisfies
\begin{align}\label{eq:Pohozaev}
&  \frac 12 \int_{\Gamma_t} |\nabla v|^2 dS=\int_{\Gamma_t}
  \left|\frac{\partial v}{\partial s}\right|^2 dS- \int_{\mathcal C_t}
  e^{-2s}\, \widetilde hv\frac{\partial v}{\partial s} \, d\mu+\frac{N-2}2
  \int_{\mathcal C_t} e^{-2s} \widetilde f(s,\theta,v)v\, d\mu\\
  \notag & -\int_{\mathcal C_t} e^{-(N+1)s} \nabla_x F(e^{-s}\theta,e^{\frac{N-2}2 s}v(s,\theta))\cdot \theta \, d\mu
  -N\int_{\mathcal C_t} e^{-Ns} F(e^{-s}\theta,e^{\frac{N-2}2 s}v(s,\theta))\, d\mu\\
  \notag &
  +\int_{\Gamma_t} e^{-Nt} F(e^{-t}\theta,e^{\frac{N-2}2 t}v(t,\theta)) \, dS
   .
\end{align}
\end{Proposition}

\begin{pf} Since $v\in H_\mu(\mathcal C_t)$ for any $t>-\log R$ then
\begin{equation*}
\int_t^{+\infty}\left(\int_{\SN} |\nabla_{\mathcal C} v(s,\theta)|^2 dS(\theta)\right) \, ds
=\int_{\mathcal C_t} |\nabla_{\mathcal C} v|^2 d\mu<+\infty
\end{equation*}
which implies that the map $s\mapsto \int_{\SN} |\nabla_{\mathcal C} v(s,\theta)|^2 dS(\theta)$ is integrable in
$(t,+\infty)$ and hence
$$
\liminf_{s\to +\infty} \int_{\SN} |\nabla_{\mathcal C}
v(s,\theta)|^2 dS(\theta)=0 \, .
$$
By Lemmas \ref{l:hardy_bound}, \ref{l:sob_cil} we also have 
\begin{align*}
&\int_t^{+\infty}\left(\int_{\SN} e^{-2s} v^2(s,\theta) dS(\theta) \right) ds+
\int_t^{+\infty}\left(\int_{\SN} e^{\left(-N+\frac{N-2}2 p\right)s} |v(s,\theta)|^p dS(\theta) \right) ds\\
&
=\int_{\mathcal C_t} e^{-2s} v^2(s,\theta) \, d\mu+ 
 \int_{\mathcal C_t} e^{\left(-N+\frac{N-2}2 p\right)s} |v|^p d\mu<+\infty 
\end{align*}
so that the maps $s\mapsto \int_{\SN} e^{-2s} v^2(s,\theta) dS(\theta)$ and $s\mapsto \int_{\SN} e^{\left(-N+\frac{N-2}2 p\right)s} |v(s,\theta)|^p dS(\theta)$
are integrable in $(t,+\infty)$ and hence
\begin{equation*}
\liminf_{s\to +\infty} \left(\int_{\SN} e^{-2s} v^2(s,\theta) dS(\theta)+\int_{\SN} e^{\left(-N+\frac{N-2}2 p\right)s} |v(s,\theta)|^p dS(\theta)   \right)
=0.
\end{equation*} 
Let $\{s_k\}\subset \R$ be an increasing sequence such that $s_k\to+\infty$ and
\begin{equation} \label{eq:s_k}
\lim_{k\to+\infty} \int_{\SN} \left(|\nabla_{\mathcal C} v(s_k,\theta)|^2 
+e^{-2s_k} v^2(s_k,\theta)
+e^{\left(-N+\frac{N-2}2 p\right)s_k} |v(s_k,\theta)|^p
\right) dS(\theta)=0 \, .
\end{equation}
From Remark \ref{rem:reg}, $u\in C^1(\Omega\setminus\{0\})$ and hence $v\in C^1(\mathcal C_\Omega)$.
Since $v$ is a weak solution of \eqref{eq:equation_in_C} in $\mathcal C_\omega$ for any bounded domain $\omega$
satisfying $\partial \omega\in C^1$ and $\overline\omega\subset \Omega$, testing \eqref{eq:equation_in_C}
with $\frac{\partial v}{\partial s}$ (we recall that $\frac{\partial v}{\partial s}\in
H^1_\mu(\mathcal C_t\setminus \overline{\mathcal C_{s_k}})$ in view of
Remark \ref{rem:reg}) and using \eqref{eq:variational_Hmu} we obtain
\begin{align}\label{eq:P}
& \int_{\mathcal C_t\setminus \overline{\mathcal C}_{s_k}} \!\!\! e^{-2s}
  \, \big(\widetilde hv + \widetilde f(s,\theta,v)
\big)\frac{\partial v}{\partial s} \, d\mu  =\int_{\mathcal C_t\setminus \overline{\mathcal C}_{s_k}}\!\!\!
  \nabla_{\mathcal C} v\cdot \nabla_{\mathcal C} \left(\tfrac{\partial
      v}{\partial s}\right)  d\mu +\int_{\Gamma_t}
  \left|\frac{\partial v}{\partial s}\right|^2\!\! dS
  -\int_{\Gamma_{s_k}} \left|\frac{\partial v}{\partial s}\right|^2 \!\!dS \\
  \notag & =\frac 12 \int_{\mathcal C_t\setminus \overline{\mathcal
      C}_{s_k}} \frac{\partial}{\partial s} (|\nabla_{\mathcal C}
  v|^2) \, d\mu +\int_{\Gamma_t} \left|\frac{\partial v}{\partial
      s}\right|^2 dS
  -\int_{\Gamma_{s_k}} \left|\frac{\partial v}{\partial s}\right|^2 dS \\
  \notag & =\frac 12 \int_{t}^{s_k} \left(\frac{\partial}{\partial s}
    \int_{\SN} |\nabla_{\mathcal C} v(s,\theta)|^2 dS(\theta)\right)
  ds +\int_{\Gamma_t} \left|\frac{\partial v}{\partial s}\right|^2 dS
  -\int_{\Gamma_{s_k}} \left|\frac{\partial v}{\partial s}\right|^2 dS \\
  \notag & =\frac 12 \int_{\SN} |\nabla_{\mathcal C} v(s_k,\theta)|^2
  dS(\theta) -\frac 12 \int_{\SN} |\nabla_{\mathcal C} v(t,\theta)|^2
  dS(\theta) +\int_{\Gamma_t} \left|\frac{\partial v}{\partial
      s}\right|^2 dS -\int_{\Gamma_{s_k}} \left|\frac{\partial
      v}{\partial s}\right|^2 dS \, .
\end{align}
By \eqref{eq:s_k} we infer that
\begin{align}\label{eq:s_k_bis}
\lim_{k\to +\infty} \int_{\Gamma_{s_k}} \left|\frac{\partial v}{\partial s}\right|^2\!\! dS
&=\lim_{k\to +\infty} \int_{\SN} \left|\frac{\partial v}{\partial s}(s_k,\theta)\right|^2\!\! dS(\theta)\leq \lim_{k\to +\infty} \int_{\SN} |\nabla_{\mathcal C} v(s_k,\theta)|^2 dS(\theta)=0 .
\end{align}
Moreover an integration by parts in the left hand side of \eqref{eq:P} yields
\begin{align*}
&
\int_{\mathcal C_t\setminus \overline{\mathcal C}_{s_k}}  e^{-2s}
  \, \widetilde f(s,\theta,v)
\frac{\partial v}{\partial s} \, d\mu
=-\frac{N-2}2\int_{\mathcal C_t\setminus \overline{\mathcal C}_{s_k}} 
e^{-Ns} f(e^{-s}\theta,u(e^{-s}\theta))u(e^{-s},\theta)\, d\mu\\
&
+\int_{\mathcal C_t\setminus \overline{\mathcal C}_{s_k}} e^{-(N+1)s} 
\nabla_x F(e^{-s}\theta,u(e^{-s}\theta))\cdot \theta\, d\mu
+N\int_{\mathcal C_t\setminus \overline{\mathcal C}_{s_k}} 
e^{-Ns} F(e^{-s}\theta,u(e^{-s}\theta))\, d\mu\\
&
-\int_{\Gamma_t} e^{-Ns} F(e^{-s}\theta,u(e^{-s}\theta)) \, dS
+\int_{\Gamma_{s_k}} e^{-Ns} F(e^{-s}\theta,u(e^{-s}\theta)) \, dS \, .
\end{align*}
By \eqref{F} and \eqref{eq:s_k} we have 
\begin{align} \label{eq:s_k_ter}
\lim_{k\to +\infty} & \left|\int_{\Gamma_{s_k}} e^{-Ns} F(e^{-s}\theta,u(e^{-s}\theta)) \, dS\right|\\
\notag\le & {\rm const}\lim_{k\to +\infty} \left( \int_{\SN} e^{-2s_k} v^2(s_k,\theta)\, dS(\theta)
+
\int_{\SN} e^{\left(-N+\frac{N-2}2 p\right)s_k} |v(s_k,\theta)|^p dS(\theta)\right)=0.
\end{align}
Passing to the limit as $k\to +\infty$ in \eqref{eq:P}, by (\ref{eq:h}), (\ref{F}), \eqref{eq:s_k}, \eqref{eq:s_k_bis} and \eqref{eq:s_k_ter}
we arrive to the conclusion.
\end{pf}

\noindent In the next lemma we provide a useful representation for the derivative of $D$.
\begin{Lemma}\label{l:Dprime}
 Under the same assumptions of Proposition
  \ref{p:Pohozaev}, the function $D$ defined  in (\ref{eq:D(t)})
  belongs to $W^{1,1}_{{\rm
       loc}}(-\log R,+\infty)$
and
\begin{align*}
D'(t)=&
 -2\int_{\Gamma_t} \left|\frac{\partial v}{\partial s}\right|^2 dS
+2\int_{\mathcal C_t} e^{-2s} \widetilde h v \frac{\partial v}{\partial s} \, d\mu 
-(N-2) \int_{\mathcal C_t} e^{-2s} \widetilde f(s,\theta,v)v\, d\mu \\
&
 +2\int_{\mathcal C_t} e^{-(N+1)s} \nabla_x F(e^{-s}\theta,e^{\frac{N-2}2 s}v(s,\theta))\cdot \theta \, d\mu+2N\int_{\mathcal C_t} e^{-Ns} F(e^{-s}\theta,e^{\frac{N-2}2 s}v(s,\theta))\, d\mu\\
 &
 -2\int_{\Gamma_t} e^{-Nt} F(e^{-t}\theta,e^{\frac{N-2}2 t}v(t,\theta)) \, dS 
+ e^{-2t}\int_{\Gamma_t}\big( \widetilde h v^2 + \widetilde f(t,\theta,v)v\big) dS
\end{align*}
in a distributional sense and for a.e. $t\in (-\log R,+\infty)$.
\end{Lemma}
\begin{pf} 
Since
\begin{align*}
D'(t)&= -\int_{\Gamma_t} |\nabla_{\mathcal C} v|^2
dS+e^{-2t}\int_{\Gamma_t}\big( \widetilde h v^2+\widetilde f(t,\theta,v) v\big) dS,
\end{align*}
the proof directly follows from \eqref{eq:Pohozaev}.
\end{pf}

\noindent The derivative of $H$ is computed in the next lemma.

\begin{Lemma}\label{l:Hprime} Under the same assumptions of Proposition \ref{p:Pohozaev} let $H$ be as in \eqref{H(t)}.
Then $H$ is differentiable in $(-\log R,+\infty)$ and
\begin{equation*}
H'(t)=2 \int_{\Gamma_t} v\frac{\partial v}{\partial s}\, dS=-2D(t)
\end{equation*}
for any $t\in (\log R,+\infty)$.
\end{Lemma}

\begin{pf}
By Remark \ref{rem:reg} $v\in C^1(\mathcal C_\Omega)$. Moreover
$H(t)=\int_{\SN} v^2(t,\theta)\, dS(\theta) $
and hence
$$
H'(t)=\int_{\SN} 2 v(t,\theta) \frac{\partial v}{\partial t} \, dS(\theta)
=\int_{\Gamma_t} 2 v \frac{\partial v}{\partial s} \, dS,
$$
which, together with the identity
$$
\int_{\mathcal C_t} |\nabla_{\mathcal C} v|^2 d\mu
+\int_{\Gamma_t}  v \frac{\partial v}{\partial s} \, dS
=\int_{\mathcal C_t} e^{-2s} \big(\widetilde
hv^2+\widetilde f(s,\theta,v)v\big)
 d\mu
$$
obtained by taking $w=v$ in \eqref{eq:variational_Hmu}, completes the proof of the lemma.
\end{pf}

\noindent Let us now compute the derivative of $\mathcal N$.

 \begin{Lemma}\label{mono}
Let $h$ satisfy \eqref{eq:h}, $f$ satisfy \eqref{F},  $u\in
\H_{{\rm loc}}(\Omega)$ be a nontrivial solution of
(\ref{eq:variational}). Let
 $\mathcal N$ be the Almgren-type function defined in \eqref{N(t)}.
 Then $\mathcal N \in W^{1,1}_{{\rm
       loc}}(\bar t,+\infty)$ and
\begin{align}\label{formulona}
{\mathcal N}'(t)=\nu_1(t)+\nu_2(t)
\end{align}
in a distributional sense and for a.e. $t\in (\bar t,+\infty)$,
where
\begin{align*}
\nu_1(t):=-2 \frac{\left(\int_{\Gamma_t} \left|\frac{\partial
            v}{\partial s}\right|^2 dS\right)
      \left(\int_{\Gamma_t} v^2 dS\right)-\left(\int_{\Gamma_t} v\frac{\partial v}{\partial
            s}dS\right)^{\!2}}
  {\left(\int_{\Gamma_t} v^2 dS\right)^2}
\end{align*}
and
\begin{align*}
& \nu_2(t)=\frac{2\int_{\mathcal C_t} e^{-2s} \widetilde
h(s,\theta)v(s,\theta)\frac{\partial v}{\partial s}(s,\theta) \,
d\mu+e^{-2t}\int_{\Gamma_t} \widetilde h v^2 \,
dS}{\int_{\Gamma_t} v^2 \,
dS} \\
& \qquad
+\frac{2\int_{\mathcal C_t} e^{-(N+1)s}\nabla_xF(e^{-s}\theta,e^{\frac{N-2}2s}v(s,\theta))\cdot
\theta \, d\mu}{\int_{\Gamma_t} v^2 \, dS} \\
& \qquad +\frac{2N\int_{\mathcal C_t}
e^{-Ns}F(e^{-s}\theta,e^{\frac{N-2}2s}v(s,\theta))\, d\mu-(N-2)
\int_{\mathcal C_t} e^{-2s}
\widetilde f(s,\theta,v(s,\theta)) v(s,\theta)\,
d\mu}{\int_{\Gamma_t} v^2 \, dS}\\
& \qquad +
\frac{e^{-2t}\int_{\SN} \widetilde f(t,\theta,v(t,\theta)) v(t,\theta) \, dS(\theta)-2e^{-Nt }\int_{\SN} F(e^{-t}\theta,e^{\frac{N-2}2 t}v(t,\theta))\, dS(\theta)
}{\int_{\Gamma_t} v^2 \, dS} \, .
\end{align*}
\end{Lemma}
\begin{pf}
  It follows from \eqref{N(t)}, Lemmas \ref{l:Dprime},
  \ref{l:Hprime}.
\end{pf}

In order to show that the Almgren function $\mathcal N$ admits a
finite limit as $t\to +\infty$ we need some preliminary estimates
which will be proved in the next lemmas.

\begin{Lemma} \label{l:est-1}
Under the same assumptions as in Lemma \ref{mono}, let $\mathcal N$
be as in (\ref{N(t)}) and $\overline t$ as in Lemma \ref{welld}.
Then, up to choose a larger $\overline t$,
we have
\begin{equation*}
\mathcal N(t)\ge -\overline C e^{-M t}
\end{equation*}
and
\begin{equation} \label{eq:l:est-1}
D(t)+H(t)\ge \frac 12 \left( \int_{\mathcal C_t} |\nabla_{\mathcal
C} v|^2 d\mu+\int_{\Gamma_t} v^2 dS\right)
\end{equation}
for any 
$t>\overline t$, where $\overline C$ is a constant depending only on $N,\e,p,u, h, f$ and
$M=\min\big\{\e,\frac{2N}p-N+2\big\}$.
\end{Lemma}

\begin{pf} Combining assumptions \eqref{eq:h}, \eqref{F} with
Lemma \ref{l:hardy_bound} and Lemma \ref{l:sob_cil}, we obtain that,
for $t>\bar t$ with $\bar t$ large,
\begin{align*}
  & D(t)= \int_{\mathcal C_t} |\nabla_{\mathcal C} v|^2
  d\mu-\int_{\mathcal C_t} e^{-2s} \widetilde h v^2 d\mu-\int_{
    \mathcal C_t} e^{-2s} \widetilde f(s,\theta,v)v \, d\mu
  \\
  & \ge \left(1-C_h\widetilde C_\e e^{-\e t}-C_f \widetilde C_2
    e^{-2t} -C_f C_{N,p}
    e^{\left(-\frac{2N}p+N-2\right)t}\big({\textstyle \int_{\mathcal
        C_t} e^{\left(-N+\frac{N-2}2p\right)s} |v|^p d\mu}
    \big)^{\frac{p-2}p} \right) \int_{\mathcal C_t}
  |\nabla_{\mathcal C} v|^2 d\mu \\
  & \qquad -\left(C_h\widetilde C_\e e^{-\e t}+C_f \widetilde C_2
    e^{-2t}+C_f C_{N,p}
    e^{\left(-\frac{2N}p+N-2\right)t}\big({\textstyle \int_{\mathcal
        C_t} e^{\left(-N+\frac{N-2}2p\right)s} |v|^p d\mu}
    \big)^{\frac{p-2}p}
  \right) \int_{\Gamma_t} v^2 dS\\
  & \ge -\left(C_h\widetilde C_\e e^{-\e t}+C_f \widetilde C_2
    e^{-2t}+C_f C_{N,p}
    e^{\left(-\frac{2N}p+N-2\right)t}\big({\textstyle \int_{\mathcal
        C_t} e^{\left(-N+\frac{N-2}2p\right)s} |v|^p d\mu}
    \big)^{\frac{p-2}p} \right) \int_{\Gamma_t} v^2 dS
\end{align*} for
which yields the conclusion if $\overline t$ is chosen sufficiently
large. 
\end{pf}

Next we provide an estimate on the function $\nu_2$ introduced in
Lemma \ref{mono}.
\begin{Lemma} \label{l:est-nu2}
Under the same assumptions as in Lemma \ref{mono} we have
\begin{equation*}
|\nu_2(t)|\le \overline C_1 (e^{-\alpha t}+g(t)) (\mathcal
N(t)+1)+\overline C_2 e^{-2t} \qquad \text{for any } t>\overline t
\, ,
\end{equation*}
where $\overline C_1$ and $\overline C_2$ are two positive
constant depending only on $N,h,f,u$ but independent of $t$,
$\alpha:=\min\{\e,2,\frac{2N}p-N+2\}$, and $g\in L^1(\overline
t,+\infty)$, $g\geq 0$ a.e., satisfies 
\begin{equation*}
\int_{t}^{+\infty} g(s)\, ds\le \frac{p}{p-2} \left(\int_{\mathcal
C_{\overline t}} e^{\left(-N+\frac{N-2}2 p\right)s}
|v(s,\theta)|^p \, d\mu\right)^{\frac{p-2}p}
e^{\left(-\frac{2N}p+N-2\right)t} \qquad \text{for any }
t>\overline t \, .
\end{equation*}
\end{Lemma}
\begin{pf} 
From \eqref{eq:h} and \eqref{F} it follows that 
\begin{align} \label{eq:estimate-nu2}
|\nu_2(t)|&\le \frac{C_h\int_{\mathcal C_t} e^{-\e s}
v^2(s,\theta) \, d\mu+C_h e^{-\e t} \int_{\mathcal C_t}
|\nabla_{\mathcal C} v|^2 d\mu+C_he^{-\e t} \int_{\Gamma_t} v^2
dS} {\int_{\Gamma_t} v^2 \, dS}
\\
\notag & \qquad  + \frac{3NC_f\int_{\mathcal C_t}
e^{-Ns}[u^2(e^{-s}\theta)+|u(e^{-s}\theta)|^p] \,
d\mu}{\int_{\Gamma_t} v^2 \, dS}\\
\notag &\qquad +\frac{3C_f e^{-Nt}\int_{\SN}
[u^2(e^{-t}\theta)+|u(e^{-t}\theta)|^p]\,
dS(\theta)}{\int_{\Gamma_t} v^2 \, dS}\\
\notag &=\frac{C_h\int_{\mathcal C_t} e^{-\e s} v^2(s,\theta) \,
d\mu+C_h e^{-\e t} \int_{\mathcal C_t} |\nabla_{\mathcal C} v|^2
d\mu+C_he^{-\e t} \int_{\Gamma_t} v^2 dS} {\int_{\Gamma_t} v^2 \,
dS}\\
\notag & \qquad + \frac{3NC_f\left[\int_{\mathcal C_t}
e^{-2s}v^2(s,\theta)\, d\mu+\int_{\mathcal C_t}
e^{\left(-N+\frac{N-2}2 p\right)s}|v(s,\theta)|^p \,
d\mu\right]}{\int_{\Gamma_t} v^2 \, dS}\\
\notag &\qquad +\frac{3C_f  e^{-2t}\int_{\SN} v^2(t,\theta) \,
dS(\theta)}{\int_{\Gamma_t} v^2 \, dS}+ \frac{3C_f
e^{\left(-N+\frac{N-2}2 p\right)t}\int_{\SN} |v(t,\theta)|^p\,
dS(\theta)}{\int_{\Gamma_t} v^2 \, dS} .
\end{align}
By Lemma \ref{l:sob_cil} and \eqref{eq:l:est-1} we obtain for any
$t>\overline t$
\begin{align*}
& \left(\int_{\mathcal C_t} e^{\left(-N+\frac{N-2}2 p\right)s}
|v(s,\theta)|^p d\mu \right)^{\!\!2/p} \le C_{N,p}
e^{\left(-\frac{2N}p+N-2\right)t} \left(\int_{\mathcal C_t}
|\nabla_{\mathcal C} v|^2 d\mu+\int_{\Gamma_t} v^2 dS \right)\\
& \le 2\,C_{N,p} e^{\left(-\frac{2N}p+N-2\right)t} (D(t)+H(t))
=2\,C_{N,p} \ e^{\left(-\frac{2N}p+N-2 \right)t} (\mathcal
N(t)+1) \int_{\Gamma_t} v^2 dS
\end{align*}
and hence
\begin{align} \label{eq:g(t)}
& \frac{3\,C_f e^{\left(-N+\frac{N-2}2 p\right)t}\int_{\SN}
|v(t,\theta)|^p\, dS(\theta)}{\int_{\Gamma_t} v^2 \, dS} \\
 \notag & \qquad \qquad \le 6C_f C_{N,p}  \ e^{\left(-\frac{2N}p+N-2\right)t}
 \
\frac{\int_{\Gamma_t}e^{\left(-N+\frac{N-2}2 p\right)t}
|v(t,\theta)|^p\, dS}{\left(\int_{\mathcal C_t}
e^{\left(-N+\frac{N-2}2 p\right)s} |v(s,\theta)|^p \, d\mu
\right)^{2/p}}\, (\mathcal N(t)+1) \, .
\end{align}
We also have
\begin{align} \label{eq:integrable}
0\le g(t)&:=e^{\left(-\frac{2N}p+N-2\right)t}
 \ \frac{\int_{\Gamma_t}e^{\left(-N+\frac{N-2}2 p\right)t}
|v(t,\theta)|^p\, dS}{\left(\int_{\mathcal C_t}
e^{\left(-N+\frac{N-2}2 p\right)s} |v(s,\theta)|^p \, d\mu
\right)^{2/p}}\\
\notag & =-\frac{p}{p-2}\Bigg\{ \frac{d}{dt}\left[
e^{\left(-\frac{2N}p+N-2\right)t} \ \left(\int_{\mathcal C_t}
e^{\left(-N+\frac{N-2}2 p\right)s} |v(s,\theta)|^p \, d\mu
\right)^{\frac{p-2}p}\right]\\
\notag & \qquad
-\left(-\tfrac{2N}p+N-2\right)e^{\left(-\frac{2N}p+N-2\right)t} \
\left(\int_{\mathcal C_t} e^{\left(-N+\frac{N-2}2 p\right)s}
|v(s,\theta)|^p \, d\mu \right)^{\frac{p-2}p} \Bigg\}
\\
\notag & \le -\frac{p}{p-2} \frac{d}{dt}\left[
e^{\left(-\frac{2N}p+N-2\right)t} \ \left(\int_{\mathcal C_t}
e^{\left(-N+\frac{N-2}2 p\right)s} |v(s,\theta)|^p \, d\mu
\right)^{\frac{p-2}p}\right]
\end{align}
in the distributional sense for almost every $t>\overline t$. But
the right hand side of \eqref{eq:integrable} is integrable in
$(\overline t,+\infty)$ since
\begin{align*}
\lim_{t\to +\infty} e^{\left(-\frac{2N}p+N-2\right)t} \
\left(\int_{\mathcal C_t} e^{\left(-N+\frac{N-2}2 p\right)s}
|v(s,\theta)|^p \, d\mu \right)^{\frac{p-2}p}=0
\end{align*}
and hence we also have $g\in L^1(\overline t,+\infty)$.

Combining \eqref{eq:estimate-nu2}, \eqref{eq:g(t)} with Lemma
\ref{l:hardy_bound}, Lemma \ref{l:sob_cil} and \eqref{eq:l:est-1}
we obtain
\begin{align*}
& |\nu_2(t)|\le \Bigg[6NC_f C_{N,p}\
e^{\left(-\frac{2N}p+N-2\right)t}\left(\int_{\mathcal C_{\overline
t}} e^{\left(-N+\frac{N-2}2p\right)s} |v(s,\theta)|^p \, d\mu
\right)^{1-\frac 2p}\\
& \qquad \qquad  +2C_h(\widetilde C_\e+1) e^{-\e t}+6NC_f C_{N,2}
e^{-2t}+6C_f C_{N,p} \ g(t) \Bigg] (\mathcal N(t)+1)+3C_f e^{-2t}
\, .
\end{align*}
The statements of the lemma follow from this last estimate and the
definition of $g$.
\end{pf}

\noindent We can now prove that the function $\mathcal N$ admits a finite
limit as $t\to +\infty$.
\begin{Lemma} \label{l:convergence}
Under the same assumptions as in Lemma \ref{mono}, the limit
$\gamma:=\lim_{t\to +\infty} \mathcal N(t)$ 
exists and is finite. Moreover $\gamma\geq0$.
\end{Lemma}

\begin{pf}
From Lemma \ref{l:est-1} we have that
\begin{equation} \label{eq:lim>=0}
\liminf_{t\to +\infty} \mathcal N(t)\ge 0 \, .
\end{equation}
On the other hand, by Lemma \ref{mono}, Schwarz inequality, and
Lemma \ref{l:est-nu2}, we have
\begin{equation} \label{eq:N'}
(\mathcal N(t)+1)'=\mathcal N'(t)=\nu_1(t)+\nu_2(t)\le \nu_2(t)\le
\overline C_1(e^{-\alpha t}+g(t))(\mathcal N(t)+1)+\overline C_2
e^{-2t}
\end{equation}
and in turn
\begin{equation*}
\frac{d}{dt} \left[e^{\overline C_1 \int_t^{+\infty} (e^{-\alpha
s}+g(s))\, ds} (\mathcal N(t)+1)\right]\le \overline C_2
e^{-2t+\overline C_1 \int_t^{+\infty} (e^{-\alpha s}+g(s))\, ds}
\, .
\end{equation*}
Since the right hand side in the above line belongs to
$L^1(\overline t,+\infty)$, after integration we deduce that
$\mathcal N$ is bounded from above and hence, by \eqref{eq:N'} and
Lemma \ref{l:est-nu2}, it follows that $\mathcal N'$ is the sum of
the nonpositive function $\nu_1$ and of the integrable function
$\nu_2$. This implies that
\begin{equation*}
\lim_{t\to +\infty} \mathcal N(t)=\mathcal N(\overline
t)+\lim_{t\to +\infty} \int_{\overline t}^{t} \mathcal N'(s)\, ds
\end{equation*}
exists and it is necessarily finite since $\mathcal N$ is bounded.
This limit is necessarily nonnegative in view of
\eqref{eq:lim>=0}.
\end{pf}

\noindent As a consequence of the convergence of $\mathcal N$, the 
following estimates on $H$ hold.

\begin{Lemma}\label{l:estimates_on_H}
Suppose that all the assumptions of Lemma \ref{mono}
are satisfied. Then there exists a constant $K_1>0$ such that
\begin{equation} \label{eq:estimate-H-above}
H(\lambda)\le K_1 e^{-2\gamma\lambda} \quad \text{for any }
\lambda>\overline t,
\end{equation}
with $\gamma=\lim_{t\to+\infty}\mathcal N(t)$ as in Lemma \ref{l:convergence}.
Moreover, for any $\sigma>0$ there exists a constant $K_2(\sigma)$
such that
\begin{equation} \label{eq:estimate-H-below}
H(\lambda)\ge K_2(\sigma) e^{-(2\gamma+\sigma)\lambda} \qquad
\text{for any } \lambda>\overline t \, .
\end{equation}
\end{Lemma}
\begin{pf} By Lemma \ref{l:Hprime} and Lemma \ref{l:convergence},
we have
\begin{align*}
\frac{H'(\lambda)}{H(\lambda)}=-2\mathcal N(\lambda)=
-2\left[\gamma-\int_\lambda^{+\infty} \mathcal N'(s)\, dx\right]
\le -2\gamma +2\int_\lambda^{+\infty} \nu_2(s)\, ds \, .
\end{align*}
By Lemma \ref{l:est-nu2} we then obtain
\begin{align*}
\frac{H'(\lambda)}{H(\lambda)} \le -2\gamma +C e^{-\alpha\lambda}
\qquad \text{for any } \lambda>\overline t
\end{align*}
where $C$ is a constant depending only on $N,h,f,u, \e,N,p$ and $\alpha$
is as in Lemma \ref{l:est-nu2}.
Estimate \eqref{eq:estimate-H-above} follows after integration in
the last inequality.

On the other hand, for any $\sigma>0$ there exists
$\lambda(\sigma)>0$ such that
\begin{equation*}
\frac{H'(\lambda)}{H(\lambda)}\ge -2\gamma-\sigma \qquad \text{for
any } \lambda>\overline \lambda(\sigma) \, .
\end{equation*}
Estimate \eqref{eq:estimate-H-below} follows after integration.
\end{pf}

\section{A blow-up argument}\label{sec:blow-up-argument}
Convergence of the frequency function $\mathcal N$ as $t\to
+\infty$ is a fundamental tool in the following blow-up argument.
Hereafter, 
we denote as $0=\mu_1<\mu_2\leq \dots \le \mu_n \le\dots$ the
eigenvalues of $-\Delta_{\SN}$ with the usual notation of
repeating them as many times as their multiplicity.
 Hence we have that $\mu _{1}=\lambda _{0}=0$ and 
\begin{align*}
\text{if }k>1\text{ and }{\textstyle{\sum_{n=0}^{\ell -1}}}m_{n}<k\leq {
\textstyle{\sum_{n=0}^{\ell }}}m_{n},\text{ then }
\mu _{k}=\lambda _{\ell }.
\end{align*}

\begin{Lemma} \label{l:blow-up}
Under the same assumptions as in Lemma \ref{mono}, let us define
the family of functions $\{w_\lambda\}_{\lambda> \overline t}$
\begin{equation*}
w_\lambda(t,\theta):=\frac{v(t+\lambda,\theta)}{\sqrt{H(\lambda)}}
\qquad \text{for any } t\ge  0 \ \text{and} \ \theta\in \SN \, .
\end{equation*}
Let $\gamma$ be the limit introduced in Lemma \ref{l:convergence}.
Then
\begin{itemize}
\item[(i)] there exists $k_0\in \N\setminus \{0\}$ such that
$\gamma=\sqrt{\mu_{\kappa_0}}$ ;

\item[(ii)] for any sequence $\lambda_n\to +\infty$ there exists a
subsequence $\lambda_{n_k}$ and an eigenfunction $\psi$ of
$-\Delta_{\SN}$ corresponding to the eigenvalue $\mu_{k_0}$ such
that $\|\psi\|_{L^2(\SN)}=1$ and $\{w_{\lambda_{n_k}}\}$ converges
to the function $w(s,\theta):= e^{-\sqrt{\mu_{k_0}}s}
\psi(\theta)$ weakly in $H_\mu(\mathcal C_0)$, strongly in $H_\mu(\mathcal C_t)$ for any $t>0$, and
strongly in $C^{1,\alpha}_{{\rm loc}}(\mathcal C_0)$ for any
$\alpha\in (0,1)$.
\end{itemize}
\end{Lemma}

\begin{pf} By \eqref{eq:l:est-1} and the definition of $w_\lambda$
we see that
\begin{align*}
\mathcal N(\lambda)+1\ge \frac 12 \left(\int_{\mathcal C_0}
|\nabla_{\mathcal C} w_\lambda|^2 d\mu+\int_{\Gamma_0} w^2_\lambda
\, dS\right)
\end{align*}
and by Lemma \ref{l:equiv_norm} and Lemma \ref{l:convergence} we
deduce that $\{w_\lambda\}_{\lambda>\overline t}$ \ is bounded in
$H_\mu(\mathcal C_0)$. Let $\{\lambda_n\}_{n\in \N}$ be a sequence
such that $\lambda_n\to +\infty$. Then there exists a subsequence
$\lambda_{n_k}$ and a function $w\in H_\mu(\mathcal C_0)$ such
that $w_{\lambda_{n_k}}\rightharpoonup w$ in $H_\mu(\mathcal
C_0)$. Let us show that $w$ weakly solves the equation
$\Delta_{\mathcal C} w=0$ in $\mathcal C_0$.

By direct computation one sees that $w_\lambda$ weakly solves the
equation
\begin{equation}\label{eq:rescaled-eq}
  -\Delta_{\mathcal C} w_\lambda(t,\theta)=e^{-2\lambda} e^{-2t}
  \widetilde h(t+\lambda,\theta)
  w_\lambda(t,\theta)+\frac{e^{-2\lambda}}{\sqrt{H(\lambda)}} \,
  e^{-2t} \widetilde
  f(t+\lambda,\theta,\sqrt{H(\lambda)}w_\lambda(t,\theta)), \quad
  \text{in } \mathcal C_0,
\end{equation}
and hence, for any $\phi\in H_{\mu,0}(\mathcal C_0)$,
\begin{align} \label{eq:identita-var}
\int_{\mathcal C_0} \nabla_{\mathcal C} w_\lambda\cdot &
\nabla_{\mathcal C} \phi \, d\mu=e^{-2\lambda} \int_{\mathcal C_0}
e^{-2t} \widetilde
h(t+\lambda,\theta)w_\lambda(t,\theta)\phi(t,\theta) \, d\mu \\
\notag & +\frac{e^{-2\lambda}}{\sqrt{H(\lambda)}} \int_{\mathcal
C_0} e^{-2t} \widetilde
f(t+\lambda,\theta,\sqrt{H(\lambda)}w_\lambda(t,\theta))
\phi(t,\theta) \, d\mu \, .
\end{align}
We estimate the two terms in the right hand side of
\eqref{eq:identita-var}. By \eqref{eq:h}, Lemma
\ref{l:hardy_bound}, and boundedness of $\{w_\lambda\}$ in
$H_\mu(\mathcal C_0)$, we have
\begin{align} \label{eq:lin-term}
& \left|e^{-2\lambda} \int_{\mathcal C_0} e^{-2t} \widetilde
h(t+\lambda,\theta)w_\lambda(t,\theta)\phi(t,\theta) \, d\mu
\right| \le C_h e^{-\e\lambda} \int_{\mathcal C_0} e^{-\e t}
|w_\lambda(t,\theta)|\, |\phi(t,\theta)| \, d\mu \\
\notag & \le C_h \widetilde C_\e e^{-\e\lambda}
\left(\int_{\mathcal C_0} |\nabla_{\mathcal C} w_\lambda|^2 \,
d\mu+\int_{\Gamma_0} w_\lambda^2 dS\right)^{1/2}\cdot
\left(\int_{\mathcal C_0} |\nabla_{\mathcal C} \phi|^2 \,
d\mu+\int_{\Gamma_0} \phi^2 dS\right)^{1/2}\to 0
\end{align}
as $\lambda\to +\infty$.
On the other hand by \eqref{F} we have
\begin{align} \label{eq:nonlin-term}
&
\left|\frac{e^{-2\lambda}}{\sqrt{H(\lambda)}} \int_{\mathcal C_0}
e^{-2t} \widetilde
f(t+\lambda,\theta,\sqrt{H(\lambda)}w_\lambda(t,\theta))
\phi(t,\theta) \, d\mu\right|
\\
\notag  & \qquad \qquad \le C_f e^{-2\lambda} \int_{\mathcal C_0}
e^{-2t} |w_\lambda(t,\theta|\, |\phi(t,\theta)| \, d\mu \\
\notag  & \qquad \qquad  \qquad +C_f e^{\frac{p(N-2)-2N}2
\lambda}\int_{\mathcal C_0} e^{\frac{p(N-2)-2N}2 \, t}
|v(t+\lambda,\theta)|^{p-2} |w_\lambda(t,\theta)|\,
|\phi(t,\theta)| \, d\mu.
\end{align}
One can show that the first term in the right hand side of
\eqref{eq:nonlin-term} tends to zero as $\lambda\to +\infty$ by
proceeding as in \eqref{eq:lin-term}. Let us prove that also the
second term in right hand side of \eqref{eq:nonlin-term} tends to
zero. Indeed, by H\"older inequality, Lemma \ref{l:sob_cil} and
boundedness of $\{w_\lambda\}$ in $H_\mu(\mathcal C_0)$ we obtain
\begin{align} \label{eq:nonlin-term-2}
& e^{\frac{p(N-2)-2N}2 \lambda}\int_{\mathcal C_0}
e^{\frac{p(N-2)-2N}2 t} |v(t+\lambda,\theta)|^{p-2}
|w_\lambda(t,\theta)|\, |\phi(t,\theta)| \, d\mu\\
\notag  & \qquad \qquad \le e^{\frac{p(N-2)-2N}2 \lambda} \left(
\int_{\mathcal C_0} e^{\frac{p(N-2)-2N}2 t}
|v(t+\lambda,\theta)|^p \, d\mu\right)^{\frac{p-2}p} \times
\\
\notag  &\qquad \qquad \qquad \times \left( \int_{\mathcal C_0}
e^{\frac{p(N-2)-2N}2 t} |w_\lambda(t,\theta)|^p \,
d\mu\right)^{\frac 1p} \left(\int_{\mathcal C_0}
e^{\frac{p(N-2)-2N}2 t} |\phi(t,\theta)|^p \, d\mu\right)^{\frac
1p}
\\
\notag & \qquad \qquad \le C_{N,p} e^{\frac{p(N-2)-2N}2 \lambda}
e^{\frac{p-2}p \frac{2N-p(N-2)}2 \lambda} \left( \int_{\mathcal
C_\lambda}
e^{\frac{p(N-2)-2N}2 t} |v|^p d\mu\right)^{\frac{p-2}p}\\
\notag  &\qquad \qquad \qquad \times \left( \int_{\mathcal C_0}
|\nabla_{\mathcal C} w_\lambda|^2 d\mu+\int_{\Gamma_0} w_\lambda^2
\, dS \right)^{\frac 12} \left( \int_{\mathcal C_0}
|\nabla_{\mathcal C} \phi|^2 d\mu+\int_{\Gamma_0} \phi^2 \, dS
\right)^{\frac 12}\to 0^+
\end{align}
as $\lambda\to +\infty$ since $p<2^*$. Passing to the limit in
\eqref{eq:identita-var} along the sequence $\{\lambda_{n_k}\}$ and
using \eqref{eq:lin-term}, \eqref{eq:nonlin-term},
\eqref{eq:nonlin-term-2} we obtain
\begin{equation*}
\int_{\mathcal C_0} \nabla_{\mathcal C} w\cdot \nabla_{\mathcal C}
\phi \, d\mu=0 \qquad \text{for any } \phi\in H_{\mu,0}(\mathcal
C_0) \, .
\end{equation*}
This proves that $w$ is harmonic in $\mathcal C_0$. Testing with a
function $\phi\in H_\mu(\mathcal C_0)$ we also have
\begin{equation} \label{eq:test-phi}
\int_{\mathcal C_0} \nabla_{\mathcal C} w\cdot \nabla_{\mathcal C}
\phi \, d\mu=-\int_{\Gamma_0} \frac{\partial w}{\partial s} \,
\phi \, dS \, .
\end{equation}
Moreover, since $\int_{\Gamma_0} w_\lambda^2 \, dS=1$,  by compactness of the trace map we also have that
\begin{equation} \label{eq:Gamma-0}
\int_{\Gamma_0} w^2 \, dS=1.
\end{equation}
We claim that $w_{\lambda_{n_k}}\to w$ strongly in
$H_\mu(\mathcal C_t)$ for any $t>0$. To this aim, we first observe
that, by direct computation,  the function
$u_\lambda=T^{-1} w_\lambda$ is actually a rescaling of the function $u$, i.e.
$$
u_{\lambda}(x)=\frac{e^{-\frac{N-2}2 \lambda}}{\sqrt{H(\lambda)}} \, u(e^{-\lambda}x)
$$
so that it solves the equation
\begin{align*}
  & -\Delta u_{\lambda}-\left(\frac{N-2}2\right)^{\!\!2}
  \frac{u_\lambda}{|x|^2}=e^{-2\lambda} h(e^{-\lambda} x) u_{\lambda}+
  \frac{e^{-\frac{N-2}2 \lambda}\cdot
    e^{-2\lambda}}{\sqrt{H(\lambda)}} \,
  f(e^{-\lambda}x,\sqrt{H(\lambda)}e^{\frac{N-2}2 \lambda} u_\lambda)
  \quad \text{in } B_1.
\end{align*}
By \eqref{eq:h} we obtain
\begin{align} \label{eq:boot-1}
|e^{-2\lambda} h(e^{-\lambda} x)u_\lambda(x)|\le C_h e^{-\e
\lambda} |x|^{-2+\e} |u_\lambda(x)|
\end{align}
and by \eqref{F} and \eqref{eq:estimate-H-above}
\begin{multline}\label{eq:boot-2}
  \left|\frac{e^{-\frac{N-2}2 \lambda}\cdot
      e^{-2\lambda}}{\sqrt{H(\lambda)}} \,
    f(e^{-\lambda}x,\sqrt{H(\lambda)}e^{\frac{N-2}2 \lambda}
    u_\lambda(x))\right|\\
  \le C_fe^{-2\lambda} |u_\lambda(x)|+C_f(H(\lambda))^{\frac{p-2}2}
  e^{\left(-N+\frac{N-2}2p\right)\lambda} |u_\lambda(x)|^{p-1} \, .
\end{multline}
Taking into account that the set $\{u_\lambda\}_{\lambda>\overline
t}$ is bounded in $\mathcal H(B_1)$ (we recall that $T$ is an
isometry), by \eqref{eq:un-um} we also have that
$\{u_\lambda\}_{\lambda>\overline t}$ is also bounded in $H^1(A)$
for any open set $A\Subset B_1\setminus\{0\}$.

Therefore by \eqref{eq:boot-1}, \eqref{eq:boot-2},
\eqref{eq:estimate-H-above}, the fact that
$p<2^*$, and a standard bootstrap argument, we deduce that
$u_\lambda$ is bounded in $C^{1,\alpha}_{{\rm loc}}(B_1\setminus
\{0\})$ for any $\alpha\in (0,1)$; the same holds true for the set
$\{w_\lambda\}$ in $C^{1,\alpha}_{{\rm loc}}(\mathcal C_0)$ for
any $\alpha\in (0,1)$.

Moreover along the subsequence $\{\lambda_{n_k}\}$ we also have
\begin{equation*}
w_{\lambda_{n_k}}\to w \qquad \text{in } C^{1,\alpha}_{{\rm
loc}}(\mathcal C_0)
\end{equation*}
for any $\alpha\in (0,1)$ and in particular
\begin{equation} \label{eq:con-C^0alpha}
\frac{\partial w_{\lambda_{n_k}}}{\partial s}\to \frac{\partial
w}{\partial s}  \qquad \text{in } C^{0,\alpha}_{{\rm
loc}}(\mathcal C_0) \, .
\end{equation}
Taking $\lambda=\lambda_{n_k}$ in \eqref{eq:rescaled-eq}, testing
in $\mathcal C_t$ with the function $w_{\lambda_{n_k}}-w\in H_\mu
(\mathcal C_t)$, for any $t>0$ we obtain
\begin{align*}
& \int_{\mathcal C_t} \nabla_{\mathcal C} w_{\lambda_{n_k}}\cdot
\nabla_{\mathcal C} (w_{\lambda_{n_k}}-w) \, d\mu=-\int_{\Gamma_t}
\frac{\partial w_{\lambda_{n_k}}}{\partial s}\,
(w_{\lambda_{n_k}}-w) \, dS \\
&+e^{-2{\lambda_{n_k}}} \int_{\mathcal C_t} e^{-2s} \widetilde
h(s+\lambda_{n_k},\theta)w_{\lambda_{n_k}}(s,\theta) (w_{\lambda_{n_k}}(s,\theta)-w(s,\theta)) \, d\mu \\
& +\frac{e^{-2\lambda_{n_k}}}{\sqrt{H(\lambda_{n_k})}}
\int_{\mathcal C_t} e^{-2s} \widetilde
f(s+\lambda_{n_k},\theta,\sqrt{H(\lambda_{n_k})}w_{\lambda_{n_k}}(s,\theta))
(w_{\lambda_{n_k}}(s,\theta)-w(s,\theta)) \, d\mu \, .
\end{align*}
Using \eqref{eq:test-phi} with $\phi=w_{\lambda_{n_k}}-w$, the
last identity then gives
\begin{align*}
& \int_{\mathcal C_t} |\nabla_{\mathcal C}(w_{\lambda_{n_k}}-w)|^2
\, d\mu=-\int_{\Gamma_t} \left(\frac{\partial
w_{\lambda_{n_k}}}{\partial s}-\frac{\partial w}{\partial
s}\right) \,
(w-w_{\lambda_{n_k}}) \, dS \\
&+e^{-2{\lambda_{n_k}}} \int_{\mathcal C_t} e^{-2s} \widetilde
h(s+\lambda_{n_k},\theta)w_{\lambda_{n_k}}(s,\theta) (w_{\lambda_{n_k}}(s,\theta)-w(s,\theta)) \, d\mu \\
& +\frac{e^{-2\lambda_{n_k}}}{\sqrt{H(\lambda_{n_k})}}
\int_{\mathcal C_t} e^{-2s} \widetilde
f(s+\lambda_{n_k},\theta,\sqrt{H(\lambda_{n_k})}w_{\lambda_{n_k}}(s,\theta))
(w_{\lambda_{n_k}}(s,\theta)-w(s,\theta)) \, d\mu \, .
\end{align*}
Passing to the limit as $k\to +\infty$, proceeding as in
\eqref{eq:lin-term}-\eqref{eq:nonlin-term-2} and using
\eqref{eq:con-C^0alpha} and the fact that $w_{\lambda_{n_k}}\to w$
in $L^2(\Gamma_t)$, we obtain $\nabla_{\mathcal C}
w_{\lambda_{n_k}}\to \nabla_{\mathcal C} w$ in $L^2(\mathcal C_t)$
and in turn, thanks to Lemma \ref{l:equiv_norm},
$w_{\lambda_{n_k}}\to w$ strongly in $H_\mu (\mathcal C_t)$ for all $t>0$.

According with \eqref{N(t)}, it is reasonable to associate to every
solution $w_\lambda$ of \eqref{eq:rescaled-eq} the Almgren-type
frequency function
\begin{equation*}
\mathcal N_\lambda(t):=\frac{D_\lambda(t)}{H_\lambda(t)} \quad
\text{for any } t\ge 0 \ \text{and} \ \lambda>\overline t,
\end{equation*}
where
\begin{align*}
D_\lambda(t):=\int_{\mathcal C_t} & |\nabla_{\mathcal C}
w_\lambda|^2 \, d \mu -e^{-2\lambda}\int_{\mathcal C_t} e^{-2s}
\widetilde
h(s+\lambda,\theta)w_\lambda^2(s,\theta)\, d\mu \\
& -\frac{e^{-2\lambda}}{\sqrt{H(\lambda)}}\int_{\mathcal C_t}
e^{-2s} \widetilde
f(s+\lambda,\theta,\sqrt{H(\lambda)}w_\lambda(s,\theta))w_\lambda(s,\theta)\,
d\mu
\end{align*}
and
\begin{equation*}
H_\lambda(t):=\int_{\Gamma_t} w_\lambda^2 dS \, .
\end{equation*}
By direct computation it follows
true:
\begin{align} \label{eq:N-trasl}
\mathcal N(t+\lambda)=\mathcal N_\lambda(t) \qquad \text{for any }
t>0 \text{ and } \lambda>\overline t  .
\end{align}
Since $w_{\lambda_{n_k}}\to w$ strongly in $H_\mu(\mathcal C_t)$
for any $t>0$, passing to the limit as $k\to \infty$ and
proceeding as in \eqref{eq:lin-term}-\eqref{eq:nonlin-term-2}, we
obtain, for any $t>0$,
\begin{equation} \label{eq:D-1}
D_{\lambda_{n_k}}(t)\to \int_{\mathcal C_t} |\nabla_{\mathcal C}
w|^2 \, d \mu
\end{equation}
and
\begin{equation} \label{eq:H-1}
H_{\lambda_{n_k}}(t)\to \int_{\Gamma_t} w^2 \, dS \, .
\end{equation}

We claim that
\begin{equation} \label{eq:norm-pos}
\int_{\mathcal C_t} |\nabla_{\mathcal C} w|^2 \, d
\mu+\int_{\Gamma_t} w^2 \, dS>0 \qquad \text{for any } t>0 \, .
\end{equation}
Indeed if there exists $t>0$ such that $\int_{\mathcal C_t}
|\nabla_{\mathcal C} w|^2 \, d \mu+\int_{\Gamma_t} w^2 \, dS=0$
then by a classical unique continuation property we deduce that
$w$ is identically zero in $\mathcal C_0$ in contradiction with
\eqref{eq:Gamma-0}.
Moreover 
\begin{equation} \label{eq:pos-Gamma}
\int_{\Gamma_t} w^2 \, dS>0
\end{equation}
 for any $t>0$ since otherwise, if there exists $t>0$ such that
$\int_{\Gamma_t} w^2 \, dS=0$, then by \eqref{eq:D-1},
\eqref{eq:H-1}, \eqref{eq:norm-pos} we would have
\begin{equation*}
\gamma=\lim_{k\to +\infty} \mathcal N(t+\lambda_{n_k})= \lim_{k\to
+\infty} \mathcal N_{\lambda_{n_k}}(t)=\lim_{k\to +\infty}
\frac{D_{\lambda_{n_k}}(t)}{H_{\lambda_{n_k}}(t)}=+\infty,
\end{equation*}
a contradiction. Therefore
\begin{align*}
\mathcal N_{\lambda_{n_k}}(t)\to \mathcal
N_{w}(t):=\frac{\int_{\mathcal C_t} |\nabla_{\mathcal C} w|^2 \, d
\mu}{\int_{\Gamma_t} w^2 dS}
\end{align*}
for any $t>0$.
Combining this with \eqref{eq:N-trasl} and Lemma
\ref{l:convergence} we deduce that
\begin{align} \label{eq:math-w}
\mathcal N_w(t)=\gamma \quad \text{for any } t>0 .
\end{align}
This means that $\mathcal N_w$ is constant and in particular, for
almost every $t>0$, by Lemma \ref{mono} we have 
\begin{align*}
0=\mathcal N'_w(t)=-2 \frac{\left(\int_{\Gamma_t} \left|\frac{\partial
            w}{\partial s}\right|^2 dS\right)
      \left(\int_{\Gamma_t} w^2 dS\right)-\left(\int_{\Gamma_t} w\frac{\partial w}{\partial
            s}dS\right)^{\!2}}
  {\left(\int_{\Gamma_t} w^2 dS\right)^2} .
\end{align*}
The condition $\left(\int_{\Gamma_t} w\frac{\partial w}{\partial s}dS\right)^{\!2}
=\left(\int_{\Gamma_t} \left|\frac{\partial w}{\partial s}\right|^2 dS\right) \left(\int_{\Gamma_t} w^2 dS\right)$
implies that, for almost every $t>0$, the functions $\theta\mapsto w(t,\theta)$
and $\theta\mapsto \frac{\partial w}{\partial s}(t,\theta)$ are parallel as vectors of $L^2(\SN)$ and hence there exists a
function $\eta$ depending only on $t$ such that
\begin{equation} \label{eq:ODE}
\frac{\partial w}{\partial t}(t,\theta)=\eta(t)w(t,\theta) \qquad
\text{for any } t>0.
\end{equation}
Clearly the function $\eta(t)=(w(t,\theta))^{-1} \frac{\partial
w}{\partial t}(t,\theta)$ is well defined and continuous for any
$t>0$ thanks to \eqref{eq:pos-Gamma}.
After integration in \eqref{eq:ODE} we deduce that $w$ admits the
representation
\begin{equation*}
w(t,\theta)=\varphi(t)\psi(\theta ).
\end{equation*}
It is not restrictive assuming that $\int_{\SN} \psi^2 dS=1$.
Inserting the above representation of $w$ into the equation
$\Delta_{\mathcal C} w=0$, it follows that there exist $k_0\in
\N\setminus\{0\}$ and $c_1,c_2\in\R$
such that
\begin{equation*}
-\Delta_{\SN} \psi(\theta)=\mu_{k_0} \psi(\theta) \quad
\text{and} \quad \varphi(t)=c_1 e^{\sqrt{\mu_{k_0}}\, t}+c_2
e^{-\sqrt{\mu_{k_0}}\, t}.
\end{equation*}
Since $w\in H_\mu(\mathcal C_0)$ and $\int_{\Gamma_0} w^2 dS=1$,
then necessarily $c_1=0$ and $c_2=1$, so that we may write
\begin{equation}\label{eq:4}
w(t,\theta)=e^{-\sqrt{\mu_{k_0}}\, t} \psi(\theta).
\end{equation}
Finally, inserting \eqref{eq:4} into \eqref{eq:math-w}, we infer that
$\gamma=\sqrt{\mu_{k_0}}$, thus completing the proof.
\end{pf}

\noindent The next lemma  provides an upper bound for the function $v$.
\begin{Lemma} \label{l:uniform-estimate}
Suppose that all the assumptions of Lemma \ref{mono} are satisfied.
Then, up to enlarge $\overline t$, there exists a constant $C$
independent of $s$ such that
\begin{equation}\label{eq:5}
  \sup_{\Gamma_s} v^2\le C H(s) \quad \text{for any } s>\overline
  t
\end{equation}
and 
\begin{equation}\label{eq:6}
\sup_{\Gamma_s} v^2\le CK_1 e^{-2\gamma s} \quad \text{for any }
s>\overline t .
\end{equation}
\end{Lemma}
\begin{pf} 
Estimate \eqref{eq:6} follows from \eqref{eq:5} and
\eqref{eq:estimate-H-above}. In order to prove \eqref{eq:5} we proceed
by contradiction and assume that there
exists a sequence $s_n\to +\infty$ such that
\begin{equation*}
\sup_{\theta\in \SN} v^2(s_n,\theta)>n \int_{\SN} v^2(s_n,\theta)
\, dS(\theta) \, .
\end{equation*}
Putting $\lambda_n:=s_n-1$ and dividing both sides of the last
inequality by $\sqrt{H(\lambda_n)}$ we infer
\begin{equation} \label{eq:pointwise}
\sup_{\theta\in \SN} w_{\lambda_n}^2(1,\theta)>n \int_{\SN}
w_{\lambda_n}^2(1,\theta) \, dS(\theta)
\end{equation}
with $w_{\lambda_n}$ as in Lemma \ref{l:blow-up}.
By Lemma \ref{l:blow-up}, along a suitable subsequence
$\{\lambda_{n_k}\}$ we have
\begin{equation*}
\sup_{\theta\in \SN} w_{\lambda_{n_k}}^2(1,\theta)\to
\sup_{\theta\in\SN} e^{-2\gamma} \psi^2(\theta)
\end{equation*}
and
\begin{equation*}
\int_{\SN} w_{\lambda_{n_k}}^2(1,\theta) \, dS(\theta)\to
e^{-2\gamma} \int_{\SN} \psi^2(\theta)\, dS(\theta)=e^{-2\gamma},
\end{equation*}
hence contradicting \eqref{eq:pointwise}.
\end{pf}

\noindent We now describe the behavior of $H(t)$ as $t \to +\infty$.

\begin{Lemma} \label{l:H-2-gamma}
Suppose that all the assumptions of Lemma \ref{mono} are satisfied
and let $\gamma$ be as in Lemma~\ref{l:convergence}. Then the
limit
\begin{equation} \label{eq:H-2-gamma}
\lim_{t\to +\infty} e^{2\gamma t} H(t)
\end{equation}
exists and belongs to $(0,+\infty)$.
\end{Lemma}
\begin{pf} By Lemma \ref{l:Hprime}, Lemma \ref{l:convergence}, and direct computations we obtain
\begin{align*}
  \frac{d}{dt} (e^{2\gamma t} H(t))=2\gamma e^{2\gamma t}
  H(t)+e^{2\gamma t} H'(t)=2 e^{2\gamma t} H(t) (\gamma-\mathcal
  N(t))=2 e^{2\gamma t} H(t) \int_t^{+\infty} \mathcal N'(s)\, ds \, .
\end{align*}
Integration in $(\overline t,t)$ then yields
\begin{align} \label{eq:ex-lim-H}
e^{2\gamma t} H(t)-e^{2\gamma \overline t} H(\overline t) = 2\!
\int_{\overline t}^t\!\! e^{2\gamma s} H(s) \left(\int_s^{+\infty}
\!\!\nu_1(z)  dz\right) ds + 2\! \int_{\overline t}^t\! e^{2\gamma s}
H(s) \left(\int_s^{+\infty}\!\! \nu_2(z) \, dz\right) ds .
\end{align}
By Lemma \ref{mono} we deduce that the function
\begin{equation*}
s\mapsto e^{2\gamma s} H(s)\int_s^{+\infty} \nu_1(z) \, dz
\end{equation*}
is non positive.

On the other hand combining Lemma \ref{l:est-nu2} with
\eqref{eq:estimate-H-above} we infer that
\begin{equation*}
s\mapsto e^{2\gamma s} H(s)\int_s^{+\infty} \nu_2(z) \, dz
\end{equation*}
is integrable in a neighborhood of infinity. This implies that the
right hand side of \eqref{eq:ex-lim-H} admits a limit as $t\to
+\infty$. This proves that the limit in \eqref{eq:H-2-gamma} exists;
on the other hand by \eqref{eq:estimate-H-above} it is necessarily
finite.  It remains to prove that it is strictly positive.

Let $R$ be such that $\overline{B}_R\subset\Omega$ and let $T:=-\log R$.
For any $k\in \N$ let us denote by $\psi_k$ an eigenfunction of
$-\Delta_{\SN}$ corresponding to the eigenvalue $\mu_k$ and
suppose that the set $\{\psi_k\}_{k\ge 1}$ is an orthonormal basis
of $L^2(\SN)$.
For any $t\ge T$, we define the functions
\begin{equation*}
\varphi_k(t):=\int_{\SN} v(t,\theta)\psi_k(\theta) \, dS(\theta)
\end{equation*}
and
\begin{equation*}
\zeta_k(t):=\int_{\SN} \left[e^{-2t} \widetilde h(t,\theta)
v(t,\theta)+e^{-2t} \widetilde f(t,\theta,v(t,\theta))\right]
\psi_k(\theta) \, dS(\theta) \, .
\end{equation*}
Since $v$ is a solution to \eqref{eq:equation_in_C} then, for any
$k\ge 1$, $\varphi_k$ solves the equation
\begin{equation*}
  -\varphi_k''(t)+\mu_k \varphi_k(t)=\zeta_k(t) \quad \text{in }
  [T,+\infty).
\end{equation*}
Integration of the above ordinary differential equation yields
\begin{align*}
\varphi_k(t)=\left(c_1^k-\int_T^t
\frac{e^{-\sqrt{\mu_k}\,s}}{2\sqrt{\mu_k}} \, \zeta_k(s)\, ds
\right) e^{\sqrt{\mu_k}\,t} +\left( c_2^k+\int_T^t
\frac{e^{\sqrt{\mu_k}\,s}}{2\sqrt{\mu_k}} \, \zeta_k(s)\,
ds\right) e^{-\sqrt{\mu_k}\,t}
\end{align*}
for some $c_1^k,c_2^k\in \R$.  Let $k_0\ge 1$ be as in Lemma
\ref{l:blow-up} so that
\begin{equation*}
\gamma:=\lim_{t\to +\infty} \mathcal N(t)=\sqrt{\mu_{k_0}}.
\end{equation*}
By definition of $\varphi_k$ and the Parseval
identity we have $H(t)=\ds{\sum_{k=1}^{+\infty}
|\varphi_k(t)|^2}$. In particular, by \eqref{eq:estimate-H-above} 
\begin{equation}\label{eq:7}
  |\varphi_k(t)|\leq \sqrt{H(t)}\leq \sqrt{K_1} e^{-\sqrt{\mu_{k_0}}\lambda}  \quad \text{for all }
t>\overline t.
\end{equation}
Let $m$ be the multiplicity of the eigenvalue $\mu_{k_0}$ and let
$j_0$ be such that
$$
\mu_{j_0}=\ldots=\mu_{k_0}=\ldots=\mu_{j_0+m-1} .
$$
Let us fix an index $i\in \{j_0,\dots,j_0+m-1\}$ and provide an
estimate for the function $\zeta_i$.
From \eqref{eq:h}, \eqref{F}, \eqref{eq:estimate-H-above}, and Lemma
\ref{l:uniform-estimate} we infer
\begin{align} \label{eq:stima-zeta-i}
 |\zeta_i&(t)| \le (C_h e^{-\e t}+C_f e^{-2t}) \sqrt{H(t)}
+C_f C^{\frac{p-1}2} \sqrt{\omega_{N-1}}
(H(t))^{\frac{p-1}2}
e^{\left(-N+\frac{N-2}2p\right)t} \\
\notag & \le \sqrt{K_1}(C_h e^{-\e t}+C_f
e^{-2t})e^{-\sqrt{\mu_{k_0}}\, t}+C_f C^{\frac{p-1}2} \sqrt{\omega_{N-1}} K_1^{\frac{p-1}2}
e^{\left(-N+\frac{N-2}2p\right)t} e^{-(p-1)\sqrt{\mu_{k_0}}\, t} .
\end{align}
Since $p>2$, the previous estimate gives
\begin{equation*}
  s\mapsto e^{\sqrt{\mu_{k_0}}\, s}\zeta_i(s)\in L^1(0,+\infty),
  \quad s\mapsto e^{-\sqrt{\mu_{k_0}}\, s}\zeta_i(s)\in
  L^1(0,+\infty).
\end{equation*}
This implies that
\begin{equation*}
\left( c_2^i+\int_T^t \frac{e^{\sqrt{\mu_{k_0}}\,
s}}{2\sqrt{\mu_{k_0}}} \, \zeta_i(s)\, ds\right)
e^{-\sqrt{\mu_{k_0}}\, t}=O(e^{-\sqrt{\mu_{k_0}}\,
t})=o(e^{\sqrt{\mu_{k_0}}\, t}) \qquad \text{as } t\to +\infty
\end{equation*}
and hence
\begin{equation*}
c_1^i-\int_T^{+\infty} \frac{e^{-\sqrt{\mu_{k_0}}\,
s}}{2\sqrt{\mu_{k_0}}} \, \zeta_i(s)\, ds=0
\end{equation*}
since otherwise we would have $\lim_{t\to+\infty}\varphi_i(t)
e^{-\sqrt{\mu_{k_0}}\, t}\neq0$, in contradiction with
\eqref{eq:7}.

Therefore we may write
\begin{align} \label{eq:fourier-1}
\varphi_{i}(t)=\left(\int_t^{+\infty}
\frac{e^{-\sqrt{\mu_{k_0}}\,s}}{2\sqrt{\mu_{k_0}}} \, \zeta_i(s)\,
ds \right) e^{\sqrt{\mu_{k_0}}\,t} +\left( c_2^i+\int_T^t
\frac{e^{\sqrt{\mu_{k_0}}\,s}}{2\sqrt{\mu_{k_0}}} \, \zeta_i(s)\,
ds\right) e^{-\sqrt{\mu_{k_0}}\,t}
\end{align}
and so by \eqref{eq:stima-zeta-i} we infer
\begin{align} \label{eq:fourier}
\varphi_{i}(t)=\left( c_2^i+\int_T^t
\frac{e^{\sqrt{\mu_{k_0}}\,s}}{2\sqrt{\mu_{k_0}}} \, \zeta_i(s)\,
ds\right)
e^{-\sqrt{\mu_{k_0}}\,t}+o(e^{-(\sqrt{\mu_{k_0}}+\delta)\,t})
\qquad \text{as } t\to +\infty
\end{align}
where $\delta=\min\{\e,2,-N+\frac{N-2}2p\}$.

Suppose by contradiction that $\ds{\lim_{t\to +\infty}
e^{2\sqrt{\mu_{k_0}}\, t} H(t)=0}$, so that, by \eqref{eq:7}, for any $k\ge 1$ we have
\begin{equation} \label{eq:o-piccolo}
\lim_{t\to +\infty} e^{\sqrt{\mu_{k_0}}\,t} \varphi_k(t)=0 \, .
\end{equation}
 Multiplying both sides of \eqref{eq:fourier} by
$e^{\sqrt{\mu_{k_0}}\,t}$ and exploiting \eqref{eq:o-piccolo} we
get
\begin{align*}
c_2^i+\int_T^{+\infty}
\frac{e^{\sqrt{\mu_{k_0}}\,s}}{2\sqrt{\mu_{k_0}}} \, \zeta_i(s)\,
ds=0
\end{align*}
and hence
\begin{align*}
\varphi_{i}(t)=-\left(\int_t^{+\infty}
\frac{e^{\sqrt{\mu_{k_0}}\,s}}{2\sqrt{\mu_{k_0}}} \, \zeta_i(s)\,
ds\right) e^{-\sqrt{\mu_{k_0}}\,t}+o(e^{-(\sqrt{\mu_{k_0}}+\delta)
\,t}) \quad \text{as } t\to +\infty  .
\end{align*}
Using again \eqref{eq:stima-zeta-i}, we finally obtain
\begin{equation} \label{eq:sop-enn}
\varphi_{i}(t)=O(e^{-(\sqrt{\mu_{k_0}}+\delta) \,t}) \quad
\text{as } t\to +\infty  .
\end{equation}
Therefore, by \eqref{eq:estimate-H-below} with $\sigma<2\delta$ and
\eqref{eq:sop-enn}, we have
\begin{align*}
\int_{\SN} w_\lambda(0,\theta) \psi_i(\theta) \, dS(\theta)
=(H(\lambda))^{-\frac 12}\varphi_i(\lambda)=o(1) \quad \text{as }
\lambda\to +\infty ,
\end{align*}
for any  $i\in \{j_0,\dots,j_0+m-1\}$. Passing to the limit as
$k\to +\infty$ along a subsequence as in Lemma \ref{l:blow-up},
then yields
$$
\int_{\SN} \psi(\theta)\psi_i(\theta) \, dS(\theta)=0 \quad
\text{for any } i\in \{j_0,\dots,j_0+m-1\}
$$
with $\psi$ as in Lemma \ref{l:blow-up}. This contradicts the fact
that $\|\psi\|_{L^2(\SN)}=1$ and that $\psi$ belongs to the space
generated by $\psi_{j_0},\dots, \psi_{j_0+m-1}$. The proof is thereby complete.
\end{pf}

\noindent We are now ready to prove the main theorem.

\begin{pfn}{Theorem \ref{t:asymptotic}}
Let $\{\psi_i\}_{i\ge 1}$
be as in the proof of Lemma \ref{l:H-2-gamma}. By Lemma
\ref{l:blow-up} and Lemma \ref{l:H-2-gamma}, for any
sequence $\lambda_n\to +\infty$ there exists a subsequence
$\{\lambda_{n_k}\}$ such that, for any $\alpha\in (0,1)$,
\begin{align} \label{eq:conv-1}
e^{\gamma \lambda_{n_k}} v(\lambda_{n_k},\theta)\to
\sum_{i=j_0}^{j_0+m-1} \beta_i \psi_i(\theta) \quad \text{in }
C^{1,\alpha}(\SN) \ \text{as } k\to +\infty  ,
\end{align}
\begin{align} \label{eq:conv-2}
e^{\gamma \lambda_{n_k}} \frac{\partial v}{\partial
t}(\lambda_{n_k},\theta)\to -\gamma  \sum_{i=j_0}^{j_0+m-1}
\beta_i \psi_i(\theta) \quad \text{in } C^{0,\alpha}(\SN) \
\text{as } k\to +\infty,
\end{align}
and
\begin{align} \label{eq:conv-3}
e^{\gamma \lambda_{n_k}} \nabla_{\SN} v(\lambda_{n_k},\theta) \to
\sum_{i=j_0}^{j_0+m-1} \beta_i \nabla_{\SN} \psi_i(\theta) \quad
\text{in } C^{0,\alpha}(\SN,T\SN) \ \text{as } k\to +\infty 
\end{align}
for some $\beta_{j_0},\dots,\beta_{j_0+m-1}\in \R$ such that
$(\beta_{j_0},\dots,\beta_{j_0+m-1})\neq(0,\dots,0)$ .

Let us prove that the coefficients
$\beta_{j_0},\dots,\beta_{j_0+m-1}\in \R$ depend neither on the
sequence $\{\lambda_n\}$ nor on its subsequence
$\{\lambda_{n_k}\}$.

First of all, for any $i\in\{j_0,\dots,j_0+m-1\}$ we have
\begin{align} \label{eq:id-beta} \lim_{k\to +\infty} & e^{\gamma
    \lambda_{n_k}} \varphi_i(\lambda_{n_k})=\lim_{k\to +\infty}
  e^{\gamma \lambda_{n_k}}\int_{\SN} v(\lambda_{n_k},\theta) \,
  \psi_i(\theta)\, dS(\theta)=\beta_i .
\end{align}
On the other hand, by \eqref{eq:fourier-1} with $t=T:=-\log R$ and $R$
as in the statement of Theorem \ref{t:asymptotic}, we infer
\begin{equation*}
c_2^i=e^{\sqrt{\mu_{k_0}}\, T}\varphi_i(T)-e^{2\sqrt{\mu_{k_0}}\,
T} \int_T^{+\infty}
\frac{e^{-\sqrt{\mu_{k_0}}\,s}}{2\sqrt{\mu_{k_0}}} \, \zeta_i(s)\,
ds,
\end{equation*}
which inserted in \eqref{eq:fourier} gives
\begin{align*}
\varphi_i(t)&=\left(e^{\sqrt{\mu_{k_0}}\,
T}\varphi_i(T)-e^{2\sqrt{\mu_{k_0}}\, T} \int_T^{+\infty}
\frac{e^{-\sqrt{\mu_{k_0}}\,s}}{2\sqrt{\mu_{k_0}}} \, \zeta_i(s)\,
ds\right) e^{-\sqrt{\mu_{k_0}}\,t}
\\
& \qquad +e^{-\sqrt{\mu_{k_0}}\,t} \int_T^t
\frac{e^{\sqrt{\mu_{k_0}}\,s}}{2\sqrt{\mu_{k_0}}} \, \zeta_i(s)\,
ds+o(e^{-(\sqrt{\mu_{k_0}}+\delta)\,t}) \qquad \text{as } t\to
+\infty \, .
\end{align*}
Multiplying both sides of the last identity by
$e^{\sqrt{\mu_{k_0}}t}$ and passing to the limit as $t\to +\infty$
we obtain
\begin{align*}
e^{\sqrt{\mu_{k_0}}\,t}\varphi_i(t)\to e^{\sqrt{\mu_{k_0}}\,
T}\varphi_i(T)-e^{2\sqrt{\mu_{k_0}}\, T} \int_T^{+\infty}
\frac{e^{-\sqrt{\mu_{k_0}}\,s}}{2\sqrt{\mu_{k_0}}} \, \zeta_i(s)\,
ds+\int_T^{+\infty}
\frac{e^{\sqrt{\mu_{k_0}}\,s}}{2\sqrt{\mu_{k_0}}} \, \zeta_i(s)\,
ds \, .
\end{align*}
This combined with \eqref{eq:id-beta} yields
\begin{equation*}
\beta_i=e^{\sqrt{\mu_{k_0}}\,
T}\varphi_i(T)-e^{2\sqrt{\mu_{k_0}}\, T} \int_T^{+\infty}
\frac{e^{-\sqrt{\mu_{k_0}}\,s}}{2\sqrt{\mu_{k_0}}} \, \zeta_i(s)\,
ds+\int_T^{+\infty}
\frac{e^{\sqrt{\mu_{k_0}}\,s}}{2\sqrt{\mu_{k_0}}} \, \zeta_i(s)\,
ds \, .
\end{equation*}
Therefore the coefficients $\beta_{j_0},\dots,\beta_{j_0+m-1}$ do
depend neither on $\{\lambda_n\}$ nor on $\{\lambda_{n_k}\}$ and hence
\eqref{eq:conv-1}--\eqref{eq:conv-3} also hold as $\lambda\to
+\infty$ and not only along the sequence $\{\lambda_{n_k}\}$.
The proof of Theorem \ref{t:asymptotic} then follows from
\eqref{eq:Tu}, the fact that $v=Tu$, and direct computations.~\end{pfn}


\begin{thebibliography}{99}


\bibitem{adim-chaud-rama} Adimurthi, N. Chaudhuri, M. Ramaswamy, {\it An improved Hardy-Sobolev
inequality and its application,} Proc. Amer. Math. Soc. 130 (2002), no. 2, 489--
505.

\bibitem{almgren} F. J. Jr. Almgren, {\it $Q$ valued functions
    minimizing Dirichlet's integral and the regularity of area
    minimizing rectifiable currents up to codimension two,}
  Bull. Amer. Math. Soc. 8 (1983), no. 2, 327--328.

\bibitem{barbi-filippa-terti}  G. Barbitas, S. Filippas, A. Tertikas, {\it Series Expansion for $L^p$ Hardy Inequalities,}
Indiana University Mathematics Journal 52 (2003), no. 1, 171--190.

\bibitem{bra-chia-cirstea-trombetti} B. Brandolini, F. Chiacchio,
F. C. C\^irstea, C. Trombetti, {\it Local behaviour of singular
solutions for nonlinear elliptic equations in divergence form,}
preprint 2012.

\bibitem{bre-dupai-tesei} H. Brezis, L. Dupaigne, A. Tesei, {\it On a semilinear elliptic
equation with inverse-square potential,}  Selecta Math. (N.S.) 11
(2005), no. 1, 1--7.

\bibitem{BK} H. Br{e}zis, T. Kato, {\it Remarks on the
    {S}chr\"odinger operator with singular complex potentials,}  J.
  Math. Pures Appl. (9) 58 (1979), no.~2, 137--151.

\bibitem{BV} H. Brezis, J. L. Vazquez, {\it Blow-up solutions of
some nonlinear elliptic problems,} Rev. Mat. Univ. Complut. Madrid
10 (1997), no. 2, 443--469.

\bibitem{CW1} F. Catrina, Z.-Q. Wang, {\it On the
    Caffarelli-Kohn-Nirenberg inequalities,}
C. R. Acad. Sci. Paris S\'er. I Math. 330 (2000), no. 6, 437--442.

\bibitem{CW2} F. Catrina, Z.-Q. Wang, {\it On the
    {C}affarelli-{K}ohn-{N}irenberg inequalities: sharp constants,
    existence (and nonexistence), and symmetry of extremal functions,}
  Comm. Pure Appl. Math. 54 (2001), no. 2, 229--258.

\bibitem{chaud} N. Chaudhuri, {\it Bounds for the best constant in an improved Hardy-Sobolev inequality,}
Z. Anal. Anwendungen 22 (2003), no. 4, 757--765.

\bibitem{cirstea} F. C. C\^irstea, {\it A complete classification of the isolated
singularities for nonlinear elliptic equations with inverse square
potentials,} preprint 2012.

\bibitem{dupaigne} L. Dupaigne, {\it A nonlinear elliptic PDE with the
    inverse square potential,} J. Anal. Math. 86 (2002), 359--398.

\bibitem{fall-musina} M. M. Fall, R. Musina {\it Sharp nonexistence
    results for a linear elliptic inequality involving Hardy and Leray
    potentials,} J. Inequal. Appl. (2011), Article ID 917201, 21
  pages, doi:10.1155/2011/917201.

\bibitem{FFT4} V. Felli, A. Ferrero, {\it Almgren-type monotonicity
    methods for the classification of behavior at corners of solutions
    to semilinear elliptic equations,} Proc. Roy. Soc. Edinburgh Sect. A, to appear.

\bibitem{FFT} V. Felli, A. Ferrero, S. Terracini, {\it Asymptotic
    behavior of solutions to Schr\"odinger equations near an isolated
    singularity of the electromagnetic potential,}
 Journal of the European Mathematical Society 13 (2011), 119--174.

\bibitem{FFT2} V. Felli, A. Ferrero, S. Terracini, {\it On the
    behavior at collisions of solutions to Schr\"odinger equations with
    many-particle and cylindrical potentials,} Discrete Contin. Dynam. Systems. 32 (2012), 3895--3956.  

\bibitem{FFT3} V. Felli, A. Ferrero, S. Terracini, {\it A note on
    local asymptotics of solutions to singular elliptic equations via
    monotonicity methods,} Milan
  J. Math (2012), doi:10.1007/s00032-012-0174-y.

\bibitem{FMT2} V. Felli, E. M. Marchini, S. Terracini, {\it On the
  behavior of solutions to {S}chr\"odinger equations with
    dipole-type potentials near the singularity,} Discrete
  Contin. Dynam. Systems 21 (2008), 91--119.

\bibitem{FG} A. Ferrero, F. Gazzola, {\it Existence of solutions for
    singular critical growth semilinear elliptic equations,}
  J. Differential Equations 177  (2001),  no. 2, 494--522.

\bibitem{garcia-peral-primo} J. Garc\'{\i}a Azorero, I. Peral, A. Primo,
{\it A borderline case in elliptic problems involving weights of
Caffarelli-Kohn-Nirenberg type,} Nonlin. Anal. 67 (2007),
1878--1894.

\bibitem{GL} N. Garofalo, F.-H.  Lin, {\it Monotonicity properties of
    variational integrals, $A\sb p$ weights and unique continuation,}
  Indiana Univ. Math. J.  35 (1986), no. 2, 245--268.

\bibitem{ghous-morad} N. Ghoussoub and A. Moradifam, {\it On the best possible remaining term in the
Hardy inequality,} Proc. Natl. Acad. Sci. USA 105 (2008), no. 37,
13746--13751.

 \bibitem{Jan} E. Jannelli,
{\it The role played by space dimension in elliptic critical
problems,}  J. Differential Equations 156 (1999), no. 2, 407--426.

\bibitem{Han-P.} P. Han, {\it Asymptotic behavior of solutions to semilinear
elliptic equations with Hardy potential,} Proc. Amer. Math. Soc.
135 (2007), no. 2, 365--372.

\bibitem{Kurata} K. Kurata, {\it A unique continuation theorem for
    the Schr\"odinger equation with singular magnetic field},
  Proc. Amer. Math. Soc., 125 (1997), no. 3, 853--860.

\bibitem{Lee} J. M. Lee, {\it Introduction to smooth manifolds,}
Graduate Texts in Mathematics 218, Springer-Verlag, New York, 2003.

\bibitem{terracini96} S. Terracini, {\it On positive entire solutions
     to a class of equations with singular coefficient and critical
     exponent,} Adv. Diff. Equa. 1 (1996), no. 2, 241--264.

 \bibitem{tesei} A. Tesei, {\it Local properties of solutions of a
     semilinear elliptic equation with inverse-square potential,}
   J. Math. Sci. (N. Y.) 149 (2008), no. 6, 1726--1740.

\bibitem{vazquez_zuazua} J. L. Vazquez, E. Zuazua, {\it The Hardy
    inequality and the asymptotic behaviour of the heat equation with
    an inverse-square potential},  J. Funct. Anal.,  173 (2000), no. 1,
    103--153.

\bibitem{wz} Z.-Q. Wang, M. Zhu, {\it {H}ardy inequalities with
    boundary terms,} Electron. J. Differential Equations 2003, No. 43.

\bibitem{wolff} T. H. Wolff, {\it A property of measures in $\R^ N$
    and an application to unique continuation,}  Geom. Funct. Anal.  2
    (1992), no. 2, 225--284.


\end{thebibliography}
\end{document}